\documentclass[reqno]{amsart}

\usepackage{tikz-cd}

\usepackage{fullpage,dsfont}

\usepackage{hyperref} 

\usepackage{parskip}
\makeatletter 
\def\thm@space@setup{%
 \thm@preskip=\parskip \thm@postskip=0pt
}
\def\th@remark{%
  \thm@headfont{\itshape}%
  \normalfont 
  \thm@preskip\parskip \thm@postskip=0pt
}
\makeatother

\usepackage[nobysame,alphabetic,initials,msc-links]{amsrefs}

\usepackage[final]{pdfpages}

\usepackage{multirow}

\usepackage{array}

\usepackage{tikz-cd}
\usetikzlibrary{knots}
\usepackage{spath3}
\usepackage{float}
\usepackage{amssymb, amsfonts, amsxtra, amsmath}
\usepackage{mathrsfs}
\usepackage{mathdots}
\usepackage{wasysym}
\usepackage[all]{xy}
\usepackage{bbm}
\usepackage{calc}
\usepackage{accents}

\theoremstyle{proclaim}
\newtheorem{Theorem}{Theorem}[section]
\newtheorem{Lem}[Theorem]{Lemma}
\newtheorem{Cor}[Theorem]{Corollary}
\newtheorem{Prop}[Theorem]{Proposition}
\theoremstyle{fancyproclaim}

\theoremstyle{statement}
\newtheorem{Rem}[Theorem]{Remark}
\newtheorem{Rems}[Theorem]{Remarks}
\newtheorem{Def}[Theorem]{Definition}

\theoremstyle{fancyproclaim}
\numberwithin{equation}{section}

\DeclareMathOperator{\Ad}{\mathrm{Ad}}

\DeclareMathOperator{\id}{\mathrm{id}}

\DeclareMathOperator{\Tr}{\mathrm{Tr}}

\DeclareMathOperator{\Imm}{\mathrm{Im}}

\DeclareMathOperator{\Ker}{\mathrm{Ker}}
\DeclareMathOperator{\Spec}{\mathrm{Spec}}

\DeclareMathOperator{\wSL}{\widetilde{\mathrm{SL}}}
\DeclareMathOperator{\Lin}{\mathrm{Lin}}

\newcommand{\cop}{\mathrm{cop}}
\newcommand{\op}{\mathrm{op}}

\newcommand{\msB}{\mathscr{B}}
\newcommand{\msC}{\mathscr{C}}
\newcommand{\msD}{\mathscr{D}}
\newcommand{\msE}{\mathscr{E}}

\newcommand{\msI}{\mathscr{I}}

\newcommand{\msL}{\mathscr{L}}

\newcommand{\mfa}{\mathfrak{a}}

\newcommand{\mfg}{\mathfrak{g}}

\newcommand{\mfk}{\mathfrak{k}}
\newcommand{\mfl}{\mathfrak{l}}
\newcommand{\mfm}{\mathfrak{m}}
\newcommand{\mfn}{\mathfrak{n}}

\newcommand{\mfsl}{\mathfrak{sl}}

\newcommand{\mfsu}{\mathfrak{su}}
\newcommand{\mft}{\mathfrak{t}}
\newcommand{\mfu}{\mathfrak{u}}

\newcommand{\mcD}{\mathcal{D}}

\newcommand{\mcH}{\mathcal{H}}

\newcommand{\mcO}{\mathcal{O}}

\newcommand{\C}{\mathbb{C}}

\newcommand{\N}{\mathbb{N}}

\newcommand{\R}{\mathbb{R}}

\newcommand{\Y}{\mathbb{Y}}
\newcommand{\Z}{\mathbb{Z}}

\newcommand{\opp}{\mathrm{op}}

\newcommand{\Span}{\mathrm{Span}}

\newcommand{\ldb}{\{}
\newcommand{\rdb}{\}}
\newcommand{\qdif}[1]{\ldb #1\rdb}

\title{Quantum $SL(2,\R)$ and its irreducible representations}
\author{K.\ De Commer and J.\ R.\ Dzokou Talla}
\address{Vrije Universiteit Brussel}
\email{kenny.de.commer@vub.be}
\email{joel.right.dzokou.talla@vub.be}

\begin{document}
\maketitle
\begin{abstract} 
We define, for $q$ a real number, a unital $*$-algebra $U_q(\mfsl(2,\R))$ quantizing the universal enveloping $*$-algebra of $\mfsl(2,\R)$. We realize this $*$-algebra $U_q(\mfsl(2,\R))$ as a $*$-subalgebra of the Drinfeld double of $U_q(\mfsu(2))$ and its dual Hopf $*$-algebra $\mcO_q(SU(2))$. More precisely, $U_q(\mfsl(2,\R))$ is generated by the coideal $*$-subalgebra $\mcO_q(K\backslash SU(2)) \subseteq \mcO_q(SU(2))$ associated to the equatorial Podle\'s sphere, and by the associated orthogonal coideal $*$-subalgebra $U_q(\mfk) \subseteq U_q(\mfsu(2))$. We then classify all the irreducible $*$-represen\-tations of $U_q(\mfsl(2,\R))$.
\end{abstract}

\section*{Introduction}

Let $\mfl$ be a real Lie algebra with complexification $\mfg = \mfl_{\C}$, so $\mfl$ can be viewed as a \emph{real form} of $\mfg$. Let $U(\mfg)$ be the universal enveloping algebra of $\mfg$. Then the real form $\mfl$ gives rise to a $*$-algebra structure on $U(\mfg)$, determined by
\begin{equation}\label{EqStarl}
X^* = -X,\qquad X\in \mfl. 
\end{equation}
We write $U(\mfl)$ for the algebra $U(\mfg)$ endowed with this particular $*$-structure.

In \cite{Twie92} it was shown that for all semisimple $\mfl$ the above real structure can be quantized: letting $U_q(\mfg)$ be the standard Drinfeld-Jimbo quantization of $U(\mfg)$ for $q$ real, there exists a Hopf $*$-algebra structure on $U_q(\mfg)$ that has $U(\mfl)$ as its classical limit. We denote $U_q^c(\mfl)$ for the Hopf algebra $U_q(\mfg)$ with the above $*$-structure. It is a quantization in the direction of the \emph{maximally compact Cartan subalgebra} of $\mfl$, in the sense that under quantization the enveloping algebra of the latter Cartan subalgebra becomes a commutative and cocommutative Hopf $*$-subalgebra of $U_q^c(\mfl)$.  

When $\mfl$ is a compact semisimple Lie algebra $\mfu$, we will always use the above Hopf $*$-algebra structure, and we then simply write $U_q(\mfu) = U_q^c(\mfu)$. This Hopf $*$-algebra $U_q(\mfu)$ can be integrated. Indeed, with $U$ the simply connected Lie group with Lie algebra $\mfu$, there exists a compact quantum group $U_q$ whose associated coordinate Hopf $*$-algebra $\mcO_q(U)$ is a $q$-deformation of the algebra $\mcO(U)$ of regular functions on $U$. The Hopf $*$-algebra $\mcO_q(U)$ is in duality with $U_q(\mfu)$, and its associated quantum group C$^*$-algebra $C^*_q(U)$ has irreducible representations in direct one-to-one correspondence with the admissible (also known as type $1$) irreducible representations of $U_q(\mfu)$. When $\mfl$ is non-compact, it is still an open problem to find a corresponding integrated version of $U_q^c(\mfl)$ in the setting of locally compact quantum groups, with concrete results so far restricted to particular cases such as the complex case and the rank 1 case, see e.g.~ \cite{KK03,VY20}. 

In this paper, we will show that a different quantization is possible for $U(\mfl)$, by quantizing in the direction of a \emph{maximally non-compact Cartan subalgebra}. The upshot is that this approach will lead to an easier setting in which to consider the integration problem (see \cite{DCDz21}), but the price we pay is that the resulting quantization $U_q(\mfl)$ will no longer be a Hopf $*$-algebra but only a coideal $*$-subalgebra. However, this does not seem to form a serious problem in practice, as the resulting $*$-algebra has by itself still an interesting representation theory. This approach is motivated substantially by the work on coideal subalgebras of $U_q(\mfg)$ associated to symmetric pairs, as developed by G. Letzter \cite{Let99} for general finite-dimensional $\mfg$ and by S. Kolb in the Kac-Moody setting \cite{Kol14}. We illustrate the above considerations in this paper for the case $\mfl=\mfsl(2,\R)$. In particular, we give a complete classification of all irreducible $*$-representations of $U_q(\mfsl(2,\R))$. We also compare the resulting theory with the classical limit as $q$ tends to $1$. Writing $\mfsl(2,\C)_{\R}$ for $\mfsl(2,\C)$ considered as a \emph{real} Lie algebra, we note that the corresponding problem for the Hopf $*$-algebra $U_q(\mfsl(2,\C)_{\R})$ has been considered in \cite{Pus93,PW94}, see also \cite{BR99}. Note that there is also a different integrable quantization of $U(\mfsl(2,\R))$ (or rather its modular double) as a locally compact quantum group in the case of $q$ unimodular \cite{ByTe03, Ip13}.

The contents of this paper are as follows. In the \emph{first section}, we develop the general theory of Drinfeld doubles of orthogonal coideals. In the \emph{second section} we introduce the $*$-algebra $U_q(\mfsl(2,\R))$. In the \emph{third section}, we classify its irreducible representations. In the \emph{fourth section}, we give an outlook to further developments. 
 
\emph{Acknowledgements}: The work of K. De Commer and J. R. Dzokou Talla was supported by the FWO grant G032919N. We thank Christian Voigt and Bob Yuncken for their comments on the preliminary version of this article. We thank the organizers of IWOTA Lancaster UK 2021, hosted by the support of EPSRC grant EP/T007524/1, for the opportunity to present this work. We also thank an anonymous referee of an earlier version of this paper for valuable comments and remarks. 

\section{Drinfeld doubles for coideals}
 
We work over the ground field $\C$. Let $A,U$ be two Hopf $*$-algebras, with coproduct, counit and antipode written as $\Delta,\varepsilon$ and $S$. We use the sumless Sweedler notation for the coproduct,
\[
\Delta(a) = a_{(1)}\otimes a_{(2)}.
\]
Assume given a \emph{unitary pairing}
\[
\tau(-,-): A \otimes U \rightarrow \C,
\]
so, writing also $\tau$ for the tensor product pairing between $A\otimes A$ and $U \otimes U$, we have for $a,b\in A$ and $h,k\in U$ that
\begin{equation}\label{EqPairStarDel}
\tau(\Delta(a),h\otimes k) = \tau(a,hk),\qquad \tau(a\otimes b,\Delta(h)) = \tau(ab,h),
\end{equation}
\begin{equation}\label{EqPairStarEps}
\tau(1,h) = \varepsilon(h),\qquad \tau(a,1) = \varepsilon(a). 
\end{equation}
\begin{equation}\label{EqPairStar}
\tau(a^*,h) = \overline{\tau(a,S(h)^*)},\qquad \tau(a,h^*) = \overline{\tau(S(a)^*,h)}.
\end{equation}
Let $I \subseteq U$ be a (unital) left coideal $*$-subalgebra,
\[
\Delta(I) \subseteq U\otimes I.
\] 
\begin{Def}
We define
\[
I^{\perp} := \{a \in A \mid \forall h\in I: \tau(a_{(1)},h)a_{(2)} = \varepsilon(h)a\}\subseteq A.
\]
\end{Def}
The following lemma is well-known in some form or other, see e.g.~ \cite[Theorem 3.1]{Let02}, \cite{LVD07,KS14}, \cite[Proposition 2.42]{DCM20}.
\begin{Lem}
The subspace $I^{\perp}  \subseteq A$ is a right coideal $*$-subalgebra. 
\end{Lem}
We call $I^{\perp}$ the coideal \emph{orthogonal} to $I$. 
\begin{proof}
It is immediate that $I^{\perp}$ is a right coideal containing $1$: if $a\in I^{\perp}$ and $\Delta(a) = \sum_i b_i \otimes a_i \in A \otimes A$ with the $a_i$ linearly independent, then for $h \in I$ we have by coassociativity of $\Delta$ that
\[
\sum_i \tau(b_{i(1)},h) b_{i(2)} \otimes a_i = \tau(a_{(1)},h)\Delta(a_{(2)}) =  \varepsilon(h) \Delta(a) = \sum_i \varepsilon(h) b_i \otimes a_i,
\]
hence $\tau(b_{i(1)},h) b_{i(2)} = \varepsilon(h)b_i$ for all $i$ by linear independance of the $a_i$. Hence each $b_i \in I^{\perp}$, so $\Delta(a)  \in I^{\perp}\otimes A$.

To see that $I^{\perp}$ is an algebra, we compute that for $a,b\in I^{\perp}$ and $h\in I$
\[
\tau((ab)_{(1)},h)(ab)_{(2)} = \tau(a_{(1)},h_{(1)})\tau(b_{(1)},h_{(2)}) a_{(2)}b_{(2)} = \varepsilon(h)ab. 
\] 
To see that $I^{\perp}$ is $*$-invariant, we first show that 
\begin{equation}\label{EqAltDefDualCoid}
a \in I^{\perp} \Leftrightarrow \forall h\in I: \tau(a_{(1)},S^{-1}(h))a_{(2)} = \varepsilon(h)a.
\end{equation}
Indeed, consider 
\[
T: A\otimes U \rightarrow U\otimes A,\qquad  a \otimes h \mapsto \tau(a_{(1)},h_{(2)})h_{(1)}\otimes a_{(2)}.
\]
As $I$ is a left coideal, it is easily seen that 
\[
a \in I^{\perp} \Leftrightarrow \forall h\in I: T(a\otimes h) = h\otimes a. 
\]
Note now that $T$ is invertible, with inverse map
\[
T^{-1}: U\otimes A \rightarrow A\otimes U,\qquad h\otimes a \mapsto \tau(a_{(1)},S^{-1}(h_{(2)}))a_{(2)}\otimes h_{(1)}. 
\]
It is then also easily verified that for $a\in A$ we have
\[
\forall h\in I: \tau(a_{(1)},S^{-1}(h))a_{(2)} = \varepsilon(h)a \Leftrightarrow \forall h\in I: T^{-1}(h\otimes a) = a\otimes h. 
\]
Putting this together, we find \eqref{EqAltDefDualCoid}. From this, the $*$-invariance of $I^{\perp}$ immediately follows: if $a\in I^{\perp}$ and $h\in I$, then
\[
\tau(a_{(1)}^*,h)a_{(2)}^*  = \overline{\tau(a_{(1)},S(h)^*)}a_{(2)}^* = (\tau(a_{(1)},S^{-1}(h^*))a_{(2)})^* = \overline{\varepsilon(h^*)}a^* = \varepsilon(h)a^*,
\]
hence $a^* \in  I^{\perp}$.
\end{proof}

Recall now the construction of the Drinfeld double $\mcD(A,U)$ of $A$ and $U$, which we will realize as the unital $*$-algebra generated by copies of $A$ and $U$ with interchange relations
\begin{equation}\label{EqCommDrinf}
\begin{aligned}
ha &= \tau(a_{(3)},h_{(1)})a_{(2)}h_{(2)}\tau(a_{(1)},S^{-1}(h_{(3)})),\\
 ah &= \tau(S^{-1}(a_{(3)}),h_{(1)})h_{(2)}a_{(2)}\tau(a_{(1)},h_{(3)}),
\end{aligned}
\end{equation}
and with coproduct inherited from $A=(A,\Delta)$ and $U^{\cop} = (U,\Delta^{\opp})$. 

\begin{Def}\label{DefDrinfDoubCoid}
Let $I$ be a left coideal $*$-subalgebra of $U$, and let $J \subseteq I^{\perp}\subseteq A$ be a right coideal $*$-subalgebra. Then we define their \emph{Drinfeld double} 
\[
\mcD(J,I)\subseteq \mcD(A,U)
\]  
to be the unital $*$-algebra in $\mcD(A,U)$ generated by $I$ and $J$.
\end{Def} 

\begin{Lem}\label{LemDrinfCoidDecomp}
The multiplication maps
\begin{equation}\label{EqDrinfCoidDecomp}
I \otimes J \rightarrow \mcD(J,I),\qquad J\otimes I \rightarrow \mcD(J,I)
\end{equation}
are bijective, and $\mcD(J,I)$ is a right coideal $*$-subalgebra of $\mcD(A,U)$. 
\end{Lem} 
\begin{proof}
It is immediate that $\mcD(J,I)$ is a right coideal $*$-subalgebra of $\mcD(A,U)$. Now for $a\in J$ and $h\in I$ the commutation relations \eqref{EqCommDrinf} simplify, using \eqref{EqAltDefDualCoid}, to
\begin{equation}\label{EqCommRelCoid}
ha = \tau(a_{(2)},h_{(1)})a_{(1)}h_{(2)},\qquad ah = \tau(S^{-1}(a_{(2)}),h_{(1)}) h_{(2)}a_{(1)}.
\end{equation}
From this, the bijectivity of the multiplications maps in \eqref{EqDrinfCoidDecomp} follows.
\end{proof} 

The relations \eqref{EqCommRelCoid} show that $\mcD(J,I)$ can also be constructed in a different way. Recall that the \emph{Heisenberg double} $\mcH(A,U)$ is the $*$-algebra generated by copies of $A$ and $U$ with multiplication law given by \eqref{EqCommRelCoid} for all $h\in U,a\in A$.  Then $\mcH(A,U)$ is a $(A\otimes U^{\cop})$-$\mcD(A,U)$-bicomodule $*$-algebra by $\Delta_A$ and $\Delta_{U}^{\op}$, and in fact a bi-Galois object \cite{Sch04}. By the above, we may also view $\mcD(J,I) \subseteq \mcH(A,U)$. 

\section{The coideal $*$-algebras $U_q(\mfsl(2,\R)_t)$}

In this section, we will use basic results on the quantized enveloping algebra of $\mfsu(2)$ and the quantized function algebra of $SU(2)$. We refer to the monograph \cite{KS97} for general information on these Hopf $*$-algebras.

\subsection{The quantized enveloping $*$-algebra $U_q(\mfsu(2))$}

Let $0<q<1$, and consider the Drinfeld-Jimbo quantized enveloping algebra $U_q(\mfsl(2,\C))$ generated by $e,f,k^{\pm 1}$ with universal relations 
\[
ke = q^2ek,\qquad kf = q^{-2}fk,\qquad ef-fe = \frac{k-k^{-1}}{q-q^{-1}}. 
\]
It becomes a Hopf algebra upon defining
\[
\Delta(e) = e\otimes 1 + k\otimes e,\qquad \Delta(f) = f\otimes k^{-1}+ 1\otimes f,\qquad \Delta(k) = k\otimes k.
\]
We can endow $U_q(\mfsl(2,\C))$ with a $*$-structure by putting 
\[
e^* = fk,\qquad f^* = k^{-1}e,\qquad k^* =k,
\]
and the resulting Hopf $*$-algebra will be denoted $U=U_q(\mfsu(2))$.

We put 
\[
\omega = ef + \frac{q^{-1}k + qk^{-1}-(q+q^{-1})}{(q^{-1}-q)^2} = fe + \frac{qk + q^{-1}k^{-1}-(q+q^{-1})}{(q^{-1}-q)^2}
\]
for the central \emph{Casimir element}. We write $(V_{n/2},\pi_{n/2})$ for the $n+1$-dimensional unitary representation of $U_q(\mfsu(2))$. Concretely, we give $V_{n/2}$ the orthonormal basis $\xi_{0},\xi_1,\ldots,\xi_n$ and let $U_q(\mfsu(2))$ act by 
\[
\pi_{n/2}(k)\xi_p =  q^{n-2p} \xi_p,
\]
\[
(q^{-1}-q)\pi_{n/2}(e)\xi_p = \sqrt{(q^{-n+p-1}-q^{n-p+1})(q^{-p}-q^p)}q^{\frac{n}{2}-p+1}\xi_{p-1},
\]
\[
(q^{-1}-q)\pi_{n/2}(f)\xi_p = \sqrt{(q^{-n+p}-q^{n-p})(q^{-p-1}-q^{p+1})} q^{-\frac{n}{2}+p} \xi_{p+1}
\]
(with $\xi_{-1} =\xi_{p+1}=0$). The Casimir then acts by the scalar 
\[
\pi_{n/2}(\omega) = \left[\frac{n}{2}\right]\left[\frac{n}{2}+1\right],
\]
where we use the standard \emph{$q$-number notation}
\begin{equation}\label{Eqqnum}
\qquad [a] = [a]_q =  \frac{q^{-a}-q^a}{q^{-1}-q}. 
\end{equation}
Upon letting $q \rightarrow 1$ in each representation, we obtain the usual representations of $U(\mfsu(2))$, interpreting $k = q^H$ with $H = \begin{pmatrix} 1 & 0 \\ 0& -1 \end{pmatrix} \in \mfsl(2,\C)$ (we will not formalize such limits explicitly in the setting of formal power series Hopf algebras, but this can easily be achieved). 

It will be convenient to use also the rescaled elements 
\[
X = q^{-1/2}(q^{-1}-q)e^*,\qquad Y = X^* = q^{-1/2} (q^{-1}-q)e,
\]
\[
\Omega = (q^{-1}-q)^2\omega +(q^{-1}+q),
\]
so that the universal relations can be rewritten in terms of $X,Y,k^{\pm 1},\Omega$ as centrality and self-adjointness of $\Omega$ and 
\[
Xk = q^2kX,\qquad Yk = q^{-2}kY,
\]
\[
XY = -1 + q \Omega k -q^2k^2,\qquad YX = -1 +q^{-1}\Omega  k -q^{-2}k^2,
\]
with 
\[
X^* = Y,\qquad k^* = k,\qquad \Omega^* = \Omega.
\]
Then 
\[
\pi_{n/2}(X)\xi_{p} = \sqrt{(1-q^{2n-2p})(q^{-2p-2}-1)}\xi_{p+1},\]
\[
\pi_{n/2}(Y)\xi_{p} = \sqrt{(1-q^{2n-2p+2})(q^{-2p}-1)}\xi_{p-1}.
\]

Fix now $t\in \R$. We recall the construction of the one-parameter family of coideals introduced in \cite{Koo93}, see also \cite{NoMi92}. Define the following elements in $U_q(\mfsu(2))$: 
\begin{equation}\label{EqStarCoidU}
B_t = q^{-1/2}(e-fk) - i(q-q^{-1})^{-1}tk,\qquad \widetilde{B}_t =  B_t + \frac{it}{q-q^{-1}}. 
\end{equation}
Then $B_t$ and $\widetilde{B}_t$ generate the same unital algebra. Since 
\begin{equation}\label{EqStarB}
B_t^* = -B_t,\qquad \widetilde{B}_t^* = -\widetilde{B}_t,
\end{equation}
this is in fact a $*$-algebra. Since $\widetilde{B}_t$ is a skew-primitive element,
\begin{equation}\label{EqSkewPrim}
\Delta(\widetilde{B}_t) = \widetilde{B}_t\otimes 1 + k \otimes \widetilde{B}_t,
\end{equation}
the polynomial algebra 
\begin{equation}\label{Eqmfkt}
U_q(\mfk_t) := \C[B_t]
\end{equation}
generated by $B_t$ is a left coideal $*$-subalgebra. Note that $\widetilde{B}_t$ allows the classical analogue
\[
b_t = E-F -\frac{it}{2}H \in \mfsu(2)
\]
where $E = \begin{pmatrix}  0 & 1 \\0 & 0 \end{pmatrix},F = \begin{pmatrix} 0 & 0 \\ 1 & 0 \end{pmatrix}$ and $H = \begin{pmatrix} 1 & 0 \\ 0 & -1\end{pmatrix}$. We view $\mfk_t = \R b_t \subseteq \mfsu(2)$. See also \cite{DCNTY20} for more on these quantizations in the formal setting. 

In the following, we will use the notation 
\begin{equation}\label{Eqqnum2}
[[a]] = [[a]]_q = q^{a}-q^{-a},\qquad \qdif{a} = \qdif{a}_q= q^{a}+q^{-a}.
\end{equation}

The next theorem is essentially given by \cite[Theorem 4.3]{Koo93}. 
\begin{Theorem}\label{TheoSpeciB}
The spectrum of $\pi_{n/2}(iB_{[[a]]})$ equals the set
\begin{equation}\label{EqFormSpecT}
\Spec(\pi_{n/2}(iB_{[[a]]})) =\left\{[a+n-2p]\mid p=0,1,\ldots,n\right\}.
\end{equation}
\end{Theorem}
\begin{proof}
Let us write $t  = [[a]]$. Then 
\[
(q-q^{-1})iB_{t} =  iX-iY +tk,
\]
so $T = (q-q^{-1})\pi_{n/2}(iB_t)$ is the selfadjoint tri-diagonal operator with
\begin{multline*}
T\xi_p = i \sqrt{(1-q^{2n-2p})(q^{-2p -2}-1)}\xi_{p+1} + tq^{n-2p} \xi_p \\ - i\sqrt{(1-q^{2n-2p+2})(q^{-2p}-1)}\xi_{p-1} 
\end{multline*}
for $0\leq p \leq n$. Introducing the new basis vectors 
\[
e_p := \frac{i^p}{\sqrt{(-1)^p(q^{-2};q^{-2})_p(q^{2n};q^{-2})_p}}\xi_p,
\]
where we use the Pochammer symbol $(a;q^{-2})_m = (1-a)(1-q^{-2}a)\ldots (1-q^{-2m+2}a)$, we obtain that 
\[
Te_p = (1-q^{2n-2p})(q^{-2p-2}-1)e_{p+1} + tq^{n-2p}e_p + e_{p-1},
\]
which now makes sense for arbitrary $p \geq 0$, i.e.\ as an operator on $\C[\N]$ or $\C^{\N}$. It then follows that inside $\C^{\N}$, we obtain for each $\lambda \in \R$ a unique normalized eigenvector $e(\lambda) = \sum_{p=0}^n q^{(a-n)p}P_p(q^{(n-a)p}\lambda)e_p$ for polynomials $P_p$ of degree $p$, normalized so that $P_0(\lambda) = 1$. The $P_p$ satisfy the following recurrence relation:
\begin{multline*}
\lambda P_p(\lambda) = P_{p+1}(\lambda) + (1-q^{-2a})q^{2n-2p}P_p(\lambda) \\ - q^{-2a+2n} (1-q^{-2p})(1-q^{-2p+2n+2})P_{p-1}(\lambda).
\end{multline*}
By e.g.\ \cite[Section 3.17]{KS98}, we see that the $P_p(\lambda)$ are dual $q^{-2}$-Krawtchouk polynomials,
\[
P_p(\lambda) = (q^{2n};q^{-2})_p K_p(\lambda;-q^{-2a},n\mid q^{-2}),\qquad p \in \{0,1,\ldots, n\}. 
\]
From the orthogonality relations given loc.~ cit. and basic theory of orthogonal polynomials, we deduce that $P_{n+1}(\lambda)$ vanishes precisely when $\lambda = q^{2p} - q^{-2a-2p+2n}$ for $p: 0 \rightarrow n$. In particular, $e(\lambda)$ has the zero coefficient at $e_{p+1}$ (and hence becomes a genuine eigenvector for the original operator $T$ on $\C^{n+1}$) precisely when 
\[
\lambda = q^{-n+a}( q^{2p} - q^{-2a-2p+2n}) = q^{2p-n+a} - q^{-2p +n-a},\qquad p: 0 \rightarrow n.
\]
By the change of variables $p \leftrightarrow n-p$, it follows that we can write the spectrum of $T =(q-q^{-1})\pi_{n/2}(iB_{[[a]]})$ as
\[
\Spec(T) = \{q^{a+n-2p}- q^{-a-n+2p}\mid p=0,1,\ldots,n\}, 
\]
leading to the desired form for the spectrum of $\pi_{n/2}(iB_{[[a]]})$. 
\end{proof}

In particular, upon writing $t = [[a]]$ (and noting that $x \mapsto [x]$ is bijective from $\R$ to $\R$), we see from \eqref{EqFormSpecT} that $t \in \Spec((q-q^{-1})\pi_{n/2}(iB_{t}))$ if and only if $n$ is even, and moreover the resulting eigenspace for $t$ is then one-dimensional. Non-zero vectors in this eigenspace are called \emph{spherical vectors}. We choose a spherical vector 
\[
\xi^{(t;n/2)} \in V_{n/2},\qquad n \in 2\Z_{\geq0},
\]
suitably normalised. In particular, in $V_1$ we have with $t = [[a]]$ that
\[
iB_{t} = \begin{pmatrix} q^2[a] & iq^{1/2}\sqrt{q+q^{-1}} & 0 \\ -iq^{1/2}\sqrt{q+q^{-1}} & [a] & iq^{-1/2}\sqrt{q+q^{-1}} \\ 0 & -iq^{-1/2}\sqrt{q+q^{-1}} & q^{-2}[a]\end{pmatrix},
\]
and we define the standard $iB_{t}$-eigenvector at $[a]$ as 
\begin{equation}\label{EqSpheric}
\xi^{(t;1)} = \begin{pmatrix} q^{-1/2} \\ (q+q^{-1})^{-1/2}it\\ q^{1/2}\end{pmatrix}
\end{equation}
The other two eigenvectors at eigenvalues $[a+2]$ and $[a-2]$ in $V_1$ are respectively 
\begin{equation}\label{EqXiPlusMin}
\xi_+^{(t;1)} =\begin{pmatrix} q^{a+1/2} \\ -i(q+q^{-1})^{1/2} \\ - q^{-a-1/2}\end{pmatrix},\qquad \xi_-^{(t;1)} =  \begin{pmatrix} -q^{-a+1/2} \\ -i(q+q^{-1})^{1/2} \\ q^{a-1/2}\end{pmatrix}.
\end{equation}

\begin{Rem}
In terms of $\widetilde{B}_t$, we can write that 
\[
\Spec(i\widetilde{B}_{[[a]]}) = \left\{\qdif{a+p}[p]\mid p\in \left\{-\frac{n}{2},-\frac{n}{2}+1,\ldots, \frac{n}{2}-1,\frac{n}{2}\right\}\right\}.
\]
\end{Rem}

\subsection{The Hopf $*$-algebra $\mcO_q(SU(2))$}\label{SubSectQuantumSU2}

Consider the Hopf $*$-algebra 
\[
A=\mcO_q(SU(2)) \subseteq \Lin(U_q(\mfsu(2)),\C)
\] 
dual to $U_q(\mfsu(2))$, generated by the entries of the spin $1/2$-representation. Hence, by definition $\mcO_q(SU(2))$ is generated by the entries of a unitary matrix
\[
U = \begin{pmatrix} \alpha &\beta\\ \gamma& \delta\end{pmatrix},
\]
which satisfies the corepresentation identity and is such that $\mcO_q(SU(2))$ is paired with $U_q(\mfsu(2))$ through the $*$-representation 
\begin{multline*}
\pi: U_q(\mfsu(2))\rightarrow M_2(\C),\qquad (U,e) = \begin{pmatrix} 0 & q^{1/2}\\ 0& 0 \end{pmatrix}, (U,f) = \begin{pmatrix} 0 & 0 \\ q^{-1/2}& 0 \end{pmatrix},\\ (U,h) = \begin{pmatrix} q & 0 \\ 0 & q^{-1}\end{pmatrix}.
\end{multline*}
It can be shown that this pairing is non-degenerate (so $\mcO_q(SU(2))$ separates the elements of $U_q(\mfsu(2))$), and that $\mcO_q(SU(2))$ is the universal algebra generated by elements $\alpha,\beta,\gamma,\delta$ such that 
\begin{subequations}\label{EqDefRelSUq}
\begin{gather}
\alpha\beta = q\beta \alpha,\quad \alpha \gamma = q\gamma \alpha,\quad \beta\gamma = \gamma \beta,\quad \beta\delta =q\delta\beta,\quad \gamma \delta = q\delta\gamma,\\
\alpha\delta - q\beta\gamma = 1,\qquad \delta \alpha -q^{-1}\gamma \beta = 1. 
\end{gather}
\end{subequations}
The $*$-structure on $\mcO_q(SU(2))$ induced by \eqref{EqPairStar} is given by 
\[
U^* = \begin{pmatrix} \alpha^* & \gamma^* \\ \beta^* & \delta^* \end{pmatrix} = \begin{pmatrix} \delta & -q^{-1}\beta \\ -q\gamma & \alpha\end{pmatrix},
\]
and the universal relations simply say that the matrix on the right is the inverse of $U$. 

Since the representations $\pi_{n/2}$ of $U_q(\mfsu(2))$ are direct summands of a tensor power of $\pi_{1/2}$, they can also be implemented by a unitary corepresentation matrix $U_{n/2}$ of $\mcO_q(SU(2))$. We write the matrix coefficients associated to the representations $\pi_{n/2}$ as $U_{n/2}$, so 
\[
(U_{n/2}(\xi,\eta),x) = \langle \xi,\pi_{n/2}(x)\eta\rangle,\qquad \xi,\eta\in \C^{n+1}. 
\]

Note that by Drinfeld duality \cite{Dri87, ES98}, we may also interpret $\mcO_q(SU(2))$ as the quantized enveloping algebra of the Lie bialgebra dual $\mfa \oplus \mfn$ of $\mfsu(2)$, where $\mfa = \R\begin{pmatrix} 1 & 0 \\ 0 & -1\end{pmatrix}$ and $\mfn = \C\begin{pmatrix} 0 & 1 \\ 0 & 0 \end{pmatrix}$. Concretely, put
 \begin{equation}\label{EqQUEGen}
\begin{aligned}
x &= q^{-1/2}\frac{\alpha -1}{q^{1/2}- q^{-1/2}},&& \quad y = \frac{\beta}{q^{1/2}-q^{-1/2}},\\ 
w &= \frac{\gamma}{q^{1/2}-q^{-1/2}},&& \quad z =q^{-1/2} \frac{\delta - 1}{q^{1/2}-q^{-1/2}}.
\end{aligned}
\end{equation}
Then the universal relations \eqref{EqDefRelSUq} become 
\[
xy - qyx  = y,\quad xw - qwx = w,\quad yw= wy, \qquad yz - qzy = y,\qquad wz - qzw = w,
\]
\[
x+z= (1-q)(xz - yw) = (1-q)(zx -q^{-2}wy)
\]
and 
\[
z^* = x,\qquad w^* = -q^{-1}y,\qquad y^* =  -q w,\qquad x^* = z.
\]
Putting $q= 1$ (and noting that we obtain then $x+z =0$, so that we can cancel $z$ from the generators), we obtain precisely the relations for $U(\mfa\oplus \mfn)$ by putting
\begin{equation}\label{EqInfVersion}
\frac{1}{2}(y+w) = \begin{pmatrix} 0 & 1 \\ 0 &0 \end{pmatrix},\quad \frac{i}{2}(y-w) = \begin{pmatrix} 0 & -i\\ 0 & 0 \end{pmatrix},\quad x = \begin{pmatrix} 1/2 & 0 \\ 0 &-1/2\end{pmatrix}.
\end{equation}
In the following, we will then switch between the notations 
\[
\mcO_q(SU(2)) = U_q(\mfa\oplus \mfn). 
\]

Let us now consider the coideal $*$-subalgebra orthogonal to $U_q(\mfk_t)$. 

\begin{Def}
Let $t\in \R$. We define $B = \mcO_q(S_t^2) = \mcO_q(K_t\backslash SU(2))$ to be the $*$-algebra generated by the matrix entries of 
\begin{equation}\label{EqDefEt}
E_t = U^*L_t U,\qquad \textrm{where }L_t = \begin{pmatrix} 0 & i \\ -i & -t\end{pmatrix}.
\end{equation}
\end{Def}
It is easily seen that $E_t$ is of the form
\[
E_t = \begin{pmatrix} q^{-1}Z_t & Y_t \\ X_t & -t-qZ_t\end{pmatrix}
\]
where $X_t,Y_t,Z_t$ satisfy 
\begin{equation}\label{DefStarQHS}
X_t^* = Y_t,\quad Z_t^*= Z_t.
\end{equation}
The $*$-algebra $\mcO_q(S_t^2)$ coincides with the function algebra of one of the Podle\'{s} quantum spheres \cite{Pod87}, and has the universal relations 
\begin{equation}\label{EqCommRelPod1}
X_tZ_t = q^2 Z_tX_t,\qquad Y_tZ_t = q^{-2}Z_tY_t,
\end{equation}
\begin{equation}\label{EqCommRelPod2} 
X_tY_t = 1-qtZ_t -q^{2}Z_t^2,\qquad Y_tX_t = 1- q^{-1}tZ_t -q^{-2}Z_t^2. 
\end{equation}
In particular, $\mcO_q(S_t^2)$ is a right coideal $*$-algebra. In fact, we have the following \cite{NoMi92}. 
\begin{Lem}
We have $\mcO_q(S_t^2) = U_q(\mfk_t)^{\perp}$. 
\end{Lem}
\begin{proof}
From the definition of $U_q(\mfk_t)^{\perp}$ and the fact that $(\id\otimes \Delta)(U) = U_{12}U_{13}$, the inclusion $\mcO_q(S_t^2) \subseteq U_q(\mfk_t)^{\perp}$ will follow once we verify that
\begin{equation}\label{EqNecComB}
(\id\otimes (-,B_t))(U^*(L_t\otimes 1)U) = \varepsilon(B_t)L_t.
\end{equation}
However, the left hand side can be simplified to 
\[
\pi_{1/2}(S(B_{t(1)})) L_t \pi_{1/2}(B_{t(2)})
\]
and \eqref{EqNecComB} follows since $U_q(\mfk_t)$ is a left coideal and $\pi_{1/2}(B_t) = -iq(q-q^{-1})^{-1}t - iL_t$ and $L_t$ commute.

Note now that the space of spherical vectors for $\pi_{n/2}$ is zero- or one-dimen\-sional according to whether $n$ is odd or even. It follows that $U_q(\mfk_t)^{\perp}$, as a representation of  $U_q(\mfsu(2))$, is equivalent to $\pi_0\oplus \pi_1\oplus \ldots$, which is also the decomposition of $\mcO_q(S_t^2)$ as a $U_q(\mfsu(2))$-representation. The reverse inclusion $\mcO_q(S_t^2) \supseteq U_q(\mfk_t)^{\perp}$ now follows.
\end{proof}

Using the notation for spherical vectors introduced at the end of the previous section, we see that $\mcO_q(S_t^2)$ has the basis
\[
\mathscr{B}_t=  \{U_{n/2}(\xi^{(t;n/2)},\xi_r)\mid n\in 2\Z_{\geq0},r=0,1,\ldots, n\}.
\]
In particular, one can verify that 
\begin{equation}\label{EqXY}
X_t = -iq^{1/2} U_1(\xi^{(t;1)},\xi_0),\qquad Y_t = iq^{-1/2} U_1(\xi^{(t;1)},\xi_2)
\end{equation}
and
\begin{equation}\label{EqZ}
 Z_t = i(q^{-1}+q)^{-1/2}U_1(\xi^{(t;1)},\xi_1) - (q^{-1}+q)^{-1}t. 
\end{equation}

Again, we will rather consider $\mcO_q(S_t^2)$ as a quantized enveloping algebra $U_q(\mfm_t)$ of some Lie subalgebra $\mfm_t \subseteq \mfa\oplus \mfn$ (see \cite{CG06} for the extension of the Drinfeld duality principle to homogeneous spaces). It is easily derived from \eqref{EqInfVersion} what this algebra corresponds to. Indeed, we compute that 
\begin{multline}\label{EqDifEt}
(q-1)^{-1}(E_t-L_t) = L_t \begin{pmatrix} x & q^{-1/2}y \\ q^{-1/2}w & z\end{pmatrix} + \begin{pmatrix} z & -q^{-3/2} y \\ -q^{1/2}w & x \end{pmatrix} L_t   \\ + (q-1) \begin{pmatrix} z & -q^{-3/2}y \\ -q^{1/2}w & x \end{pmatrix}L_t \begin{pmatrix} x & q^{-1/2}y \\ q^{-1/2}w & z\end{pmatrix},
\end{multline}
so that the classical limit corresponds to the Lie algebra generated by the entries of $\left\lbrack L_t,\begin{pmatrix} x & y\\ w & z\end{pmatrix}\right\rbrack$. This is seen to be the Lie algebra
\[
\mfm_t = \left\{\begin{pmatrix} \lambda & \gamma-it\lambda\\ 0 & -\lambda\end{pmatrix} \mid \lambda,\gamma\in \R\right\}. 
\]
As one expects, this is indeed the orthogonal complement of $\mfk_t$ under the natural real pairing of Lie bialgebras 
\[
(X,Y) = \frac{1}{2}\Imm \Tr(XY)
\]
between $\mfsu(2)$ and $\mfa\oplus \mfn$. 

In the following, we will also switch between the notations 
\[
\mcO_q(S_t^2)  = U_q(\mfm_t). 
\]

\subsection{The coideal $*$-algebra $U_q(\mfsl(2,\R)_t)$}

We apply the construction method of Definition \ref{DefDrinfDoubCoid} to the coideal $*$-algebra $U_q(\mfk_t)\subseteq U_q(\mfsu(2))$ (introduced in \eqref{Eqmfkt}) and its dual $U_q(\mfk_t)^{\perp} = U_q(\mfm_t) \subseteq U_q(\mfa\oplus \mfn)$.

\begin{Def}
We define 
\[
 U_q(\mfsl(2,\R)_t) := \mcD(U_q(\mfm_t),U_q(\mfk_t))\subseteq  \mcD(U_q(\mfa\oplus \mfn),U_q(\mfsu(2))) =: U_q(\mfsl(2,\C)_{\R}).
\]
\end{Def}

By construction, the $*$-algebra $U_q(\mfsl(2,\R)_t)$ is a right coideal $*$-subalgebra of $U_q(\mfsl(2,\C)_{\R})$. Here $U_q(\mfsl(2,\C)_{\R})$ is seen as a quantization of the universal enveloping $*$-algebra of $\mfsl(2,\C)$ as a real Lie algebra. The associated Lie subalgebra $\mfsl(2,\R)_t$ is then, according to the description in the previous section, the real Lie subalgebra of $\mfsl(2,\C)$ generated by $\mfk_t$ and $\mfm_t$, which can be realized concretely as 
\[
\mfsl(2,\R)_t = \Ad(g_t)\mfsl(2,\R),\qquad g_t = \begin{pmatrix} 1 & it/2 \\ 0 & 1 \end{pmatrix}. 
\]
In particular, $\mfsl(2,\R)_0  = \mfsl(2,\R)$. 

The decomposition 
\[
U_q(\mfsl(2,\R)_t) =  U_q(\mfk_t)U_q(\mfm_t)
\]
is an analogue of Iwasawa decomposition for $\mfsl(2,\R)$, with the solvable part grouped together.

The following lemma is immediate from Lemma \ref{LemDrinfCoidDecomp}. 
\begin{Lem}
The $*$-algebra $U_q(\mfsl(2,\R)_t)$ is the universal $*$-algebra generated by elements $X_t,Y_t,Z_t,B_t$ satisfying the relations   \eqref{EqStarB}, \eqref{DefStarQHS}, \eqref{EqCommRelPod1}, \eqref{EqCommRelPod2} and 
\begin{equation}\label{EqCommXYB}
\begin{aligned}
B_tX_t &= q^{2}X_tB_t + q((q+q^{-1})Z_t+t),\\ B_tY_t &= q^{-2} Y_tB_t +q^{-1}((q^{-1}+q)Z_t+t),
\end{aligned}
\end{equation}
\begin{equation}\label{EqCommZB}
B_tZ_t = Z_tB_t - (X_t+Y_t).
\end{equation}
\end{Lem}

We single out the following special element in $U_q(\mfsl(2,\R)_t)$.

\begin{Def}\label{DefCasisl2}
The \emph{Casimir element} of $U_q(\mfsl(2,\R)_t)$ is defined as the element
\[
\Omega_t:= iq^{-1}X_t +(q-q^{-1})iZ_tB_t-iqY_t.
\]
\end{Def}

\begin{Prop}
The Casimir element $\Omega_t$ is a central, self-adjoint element of $U_q(\mfsl(2,\R)_t)$. 
\end{Prop}
\begin{proof}
Using \eqref{EqCommZB}, we write
\begin{multline*}
\Omega_t = iq^{-1}(X_t-Z_tB_t) - iq(Y_t-Z_tB_t) = -iq^{-1}(B_tZ_t+Y_t) +iq(B_tZ_t+X_t) \\= -iq^{-1}Y_t + (q-q^{-1})iB_tZ_t+iqX_t = \Omega_t^*,
\end{multline*}
so $\Omega_t$ is self-adjoint.

It is clear on sight that $Z_t$ and $\Omega_t$ commute. Writing 
\[
-i\Omega_t = (q^{-1}X_t-qY_t) +(q-q^{-1})Z_tB_t
\]
and observing that from \eqref{EqCommXYB} 
\begin{equation}\label{EqCommBXYSum}
B_t(qY_t - q^{-1}X_t) = (q^{-1}Y_t-qX_t)B_t,
\end{equation}
it also follows from \eqref{EqCommZB} that $\Omega_t$ commutes with $B_t$. Finally, from \eqref{EqCommBXYSum}, we see that $qY_t-q^{-1}X_t$ commutes with $Z_tB_t$, and hence with $\Omega_t$. As $\Omega_t$ is self-adjoint, it also commutes with $q^{-1}Y_t - qX_t$, and so with $X_t$ and $Y_t$. It follows that $\Omega_t$ is central. 
\end{proof}

\section{Representation theory of $U_q(\mfsl(2,\R)_t)$}\label{SecRepTheoSL}

\subsection{Classification of irreducible $U_q(\mfsl(2,\R)_t)$-representations}

In this section, we classify all admissible irreducible $*$-representations of $U_q(\mfsl(2,\R)_t)$, in the sense of the following definition. We use the notation $\wSL(2,\R)$ to refer to the simply connected cover of $SL(2,\R)$. We also use again the notation as introduced in \eqref{Eqqnum} and \eqref{Eqqnum2}, and we fix $t = [[a]]$. 

\begin{Def}\label{DefIrrepSL2R}
A $U_q(\mfsl(2,\R)_t)$-module $(V,\pi)$ is called \emph{$\wSL(2,\R)$-admissible} if $\pi(iB_t)$ is diagonalizable with finite-dimensional eigenspaces. 

If $V$ is moreover a pre-Hilbert space, we say that $V$ is an \emph{$\wSL(2,\R)$-admissible representation} if the $\pi(X_t),\pi(Z_t)$ and $\pi(Y_t)$ act by bounded operators. 

We say that $V$ is an \emph{$\wSL(2,\R)$-admissible $*$-representation} if for all $\xi,\eta\in V$ and $x\in U_q(\mfsl(2,\R)_t)$ we have 
\begin{equation}\label{EqStarInv}
\langle \xi,\pi(x)\eta\rangle = \langle \pi(x^*)\xi,\eta\rangle,\qquad x\in U_q(\mfsl(2,\R)_t). 
\end{equation}

We say $V$ is an \emph{$SL(2,\R)_t$-admissible} module/representation/$*$-representa\-tion if the eigenvalues of $\pi(iB_t)$ are of the form $[c]$ for $c\in a+\Z$.
\end{Def}

\begin{Rems}
\begin{itemize}
\item The condition for $SL(2,\R)_t$-admissibility is directly modeled on Theorem \ref{TheoSpeciB}. 
\item As the $\pi(iB_t)$-eigenspaces are required to be finite-dimensional, there is no ambiguity between the notions of irreducibility in the module- and representation-theoretic setting. 
\item It follows from \eqref{EqCommRelPod2} that for any $\wSL(2,\R)$-admissible $*$-representation, the boundedness of the $\pi(X_t),\pi(Z_t)$ and $\pi(Y_t)$ is automatic. 
\end{itemize}
\end{Rems}

Our first step is based on the following observation. 

\begin{Lem}\label{LemActB}
Let $V = (V,\pi)$ be a $U_q(\mfsl(2,\R)_t)$-module. Assume $v$ is a $\pi(iB_t)$-eigenvector at eigenvalue $[c]$. Let $n \in 2\Z_{\geq0}$. Then the vector space $V(n/2)$ spanned by the $\{U_{n/2}(\xi^{(t;n/2)},\xi_p)v\mid p \in \{0,1,\ldots,n\}\}$ is invariant under $\pi(iB_t)$, and there is a map of $\C[x]$-modules 
\[
V_{n/2} \rightarrow V(n/2),\qquad \xi_p \mapsto U_{n/2}(\xi^{(n/2;t)},\xi_p)v,
\] 
where $x$ acts on the left via $\pi_{n/2}(iB_{[[c]]})$ and on the right via $\pi(iB_t)$. 
\end{Lem} 
\begin{proof}
If $z\in U_q(\mfm_t)$, we have by \eqref{EqCommRelCoid} and \eqref{EqSkewPrim} that
\begin{multline*}
(q-q^{-1})\pi(iB_tz)v = (q-q^{-1})\tau(z_{(2)},(iB_{t})_{(1)}) \pi(z_{(1)}(iB_{t})_{(2)})v \\ = i\tau(z_{(2)},X-Y)\pi(z_{(1)})v + (q-q^{-1})\tau(z_{(2)},k)\pi(z_{(1)}iB_t)v \\= (i\tau(z_{(2)},X-Y)\pi(z_{(1)}) +[[c]]\tau(z_{(2)},k)\pi(z_{(1)}))v.
\end{multline*}
Putting $\eta^{([[c]];n/2)}_p = \pi(U_{n/2}(\xi^{(t;n/2)},\xi_p))v$, we get 
\begin{multline*}
(q-q^{-1})\pi(iB_t)\eta^{([[c]];n/2)}_p = i\sqrt{(1-q^{2n-2p})(q^{-2p-2}-1)} \eta^{([[c]];n/2)}_{p+1} \\+ q^{n-2p}[[c]]\eta^{([[c]];n/2)}_{p} -i \sqrt{(1-q^{2n-2p+2})(q^{-2p}-1)} \eta^{([[c]];n/2)}_{p-1}.
\end{multline*}
So $\pi(iB_t)$ acts on the $\eta^{([[c]];n/2)}_p$ exactly as $\pi_{n/2}(iB_{[[c]]})$ acts on the $\xi_p \in V_{n/2}$.
\end{proof}

Consider then, for $c \in \R$, the following elements in $U_q(\mfm_t)$ (where we use the notation from Section \ref{SubSectQuantumSU2}, see in particular \eqref{EqXiPlusMin}):
\begin{equation}\label{EqDeftcplus}
T_c^+  := U_1(\xi^{(t;1)},\xi^{([[c]];1)}_+),
\end{equation}
\begin{equation}\label{EqDefAc}
A_{c} :=  U_1(\xi^{(t;1)},\xi^{([[c]];1)}) - (q^{-1}+q)^{-1}[[c]]t,
\end{equation}
\begin{equation}\label{EqDefTcmin}
T_c^- := (T_c^+)^* = U_1(\xi^{(t;1)},\xi^{([[c]];1)}_-),
\end{equation}
Using \eqref{EqXY} and \eqref{EqZ}, we can write 
\[
\begin{pmatrix} T_c^++t \\ A_c \\ T_c^-+t\end{pmatrix} = \begin{pmatrix}iq^{c} & -q^{-1}-q & iq^{-c} \\ iq^{-1} & q^{c}-q^{-c} & -iq \\  -iq^{-c} & -q^{-1}-q & -iq^{c}\end{pmatrix} \begin{pmatrix} X_t \\ Z_t \\Y_t\end{pmatrix}.
\]
It is clear from the definition that the span of $1,X_t,Y_t,Z_t$ is equal to the span of $1,T_c^+,A_c,T_c^-.$ Concretely, we have
\begin{equation}\label{EqXYZT}
\begin{aligned}
\qdif{c-1}\qdif{c}\qdif{c+1}&\begin{pmatrix} X_t \\ Z_t \\Y_t \end{pmatrix} =  \\  &\begin{pmatrix} -iq^{c+1} & -i(q^{-1}+q) & iq^{-c+1}\\ -1 & q^{c}-q^{-c} & -1 \\- iq^{-c-1} & i(q^{-1}+q)&iq^{c-1}\end{pmatrix}   \begin{pmatrix} \qdif{c-1}(T_c^++t)\\ \qdif{c}A_c\\ \qdif{c+1}(T_c^-+t)\end{pmatrix}.
\end{aligned}
\end{equation}

From Lemma \ref{LemActB}, we immediately obtain the following.
\begin{Prop}\label{PropActAc}
Let $V$ be an $U_q(\mfsl(2,\R)_t)$-module, and let $V_c$ be the eigenspace at eigenvalue $[c]$ for $iB_t$. Then
\[
A_cV_c \subseteq V_c,\qquad T_c^+V_{c} \subseteq V_{c+2},\qquad T_c^-V_c \subseteq V_{c-2}. 
\]
\end{Prop} 

Let us now consider the following particular types of $U_q(\mfsl(2,\R)_t)$-modules. Recall the Casimir element introduced in Definition \ref{DefCasisl2}. 

\begin{Def}
Let $\lambda \in \R$. We call an $U_q(\mfsl(2,\R)_t)$-module $V=(V,\pi)$ \emph{$\lambda$-basic} if $V$ admits a cyclic $\pi(iB_t)$-eigenvector and $\pi(\Omega_t) = \lambda \id_V$. We call a module $V$ \emph{basic} if it is $\lambda$-basic for some $\lambda$.
\end{Def}

\begin{Lem}\label{LemEigvA}
A $U_q(\mfsl(2,\R)_t)$-module $(V,\pi)$ is $\lambda$-basic if and only if there exists $b\in \R$ and a cyclic $iB_t$-eigenvector $v$ at eigenvalue $[b]$ such that $A_bv = \lambda v$. In this case, any $iB_t$-eigenvector at eigenvalue $[c]$ is a $\lambda$-eigenvector of $A_c$. 
\end{Lem}
\begin{proof}
It follows immediately from the definition of $A_c$ that $\Omega_t v = A_cv$ for $v$ an eigenvector for $iB_t$ at eigenvalue $[c]$. 
\end{proof}

To get at some structural results about basic modules, we need the following commutation relations between the $T_c^+,A_c$ and $T_c^-$. These relations are checked by direct computations (only two relations need to be computed, the other ones following by using that $\overline{A_c} = A_{-c}$ and $\overline{T_c^+} = T_{-c}^-$ with the complex conjugation leaving the generators $X_t,Z_t,Y_t$ invariant). 

\begin{Lem}\label{LemCommRelATT}
The following relations hold in $U_q(\mfm_t)$ for $t = [[a]]$:
\begin{align*}
T_{c+2}^-T_c^+ &=-A_c^2 +t(q^{-c-1}-q^{c+1})A_c + t^2 + (q^{-c-1}+q^{c+1})^2\\
&= (q^{c-a+1}+q^{-c+a-1}-A_c)(q^{c+a+1}+q^{-c-a-1}+A_c),
\end{align*}
\[
A_{c+2}T_c^+ = T_c^+A_c,\qquad A_{c-2}T_c^- = T_c^-A_c,
\]
\begin{align*}
T_{c-2}^+T_c^- &= -A_c^2 +t(q^{-c+1}-q^{c-1})A_c + t^2+ (q^{-c+1}+q^{c-1})^2\\
& = (q^{c-a-1}+q^{-c+a+1}-A_c)(q^{c+a-1}+q^{-c-a+1}+A_c).
\end{align*}
\end{Lem}

In the following, we will use the notation 
\[
T_{c}^{+,(n)} = T_{c+2n-2}^+T_{c+2n-4}^+\cdots T_{c+2}^+T_c^+,\qquad T_{c}^{-,(n)} =  T_{c-2n+2}^-T_{c-2n+4}^-\cdots T_{c-2}^-T_c^-,
\]
where by convention $T_c^{+,(0)} = T_c^{-,(0)} = 1$.

\begin{Prop}\label{PropBasicAdm}
Every basic $U_q(\mfsl(2,\R)_t)$-module is $\wSL(2,\R)$-admissible, and has one-dimensional $iB_t$-eigen-spaces. Conversely, every \emph{irreducible} $\wSL(2,\R)$-admissi\-ble $U_q(\mfsl(2,\R)_t)$-module is basic. 
\end{Prop} 
\begin{proof}
If $V$ is an irreducible $\wSL(2,\R)$-admissible module, Schur's lemma gives that the central element $\Omega_t$ must act by a scalar. As any eigenvector of $iB_t$ is then cyclic by the irreducibility condition, one direction is shown. 

Conversely, if $V$ is a $\lambda$-basic module and $v$ is a cyclic eigenvector of $iB_t$ at eigenvalue $[c]$, we claim that the $T_{c}^{\pm,(n)}v$ span $V$. Indeed, let $W$ be the span of these elements. It is sufficient to show that $W$ is an $U_q(\mfsl(2,\R)_t)$-submodule. But put $w = T_{c}^{+,(n)}v$. Then clearly $T_{c+2n}^+w \in W$. As $v$ is an $A_c$-eigenvector by Lemma \ref{LemEigvA}, we deduce by induction from Lemma \ref{LemCommRelATT} that $A_{c+2n}w\in W$ and $T_{c+2n}^-w \in W$. This implies that $\{X_tw,Z_tw,Z_tw\} \subseteq W$. A similar result holds for the $T_c^{-,(n)}v$. As $T_{c}^{\pm,(n)}v$ is an eigenvector of $iB_t$ at eigenvalue $[c\pm 2n]$ when non-zero by Proposition \ref{PropActAc}, we are done. 
\end{proof} 

There is a canonical construction of a $\lambda$-basic module as follows. Let $b\in \R$, and let $\C_{b}$ be the one-dimensional $U_q(\mfk_t)$-module given by the character 
\[
\chi_b: iB_t \mapsto [b].
\]
From \eqref{EqDrinfCoidDecomp}, it follows that we can identify the vector spaces
\begin{equation}\label{EqConcId}
U_q(\mfsl(2,\R)_t) \underset{U_q(\mfk_t)}{\otimes} \C_b \cong U_q(\mfm_t).
\end{equation}
Let $I_{\lambda,b}$ be the left ideal of $U_q(\mfm_t)$ generated by $A_b -\lambda$. From the discussion under Lemma \ref{LemActB} and Proposition \ref{PropActAc}, we see that $I_{\lambda,b}$ is an $U_q(\mfsl(2,\R)_t)$-submodule of $U_q(\mfm_t)$. Put 
\[
M_{\lambda,b} = U_q(\mfm_t)/I_{\lambda,b},
\]
and let $\widehat{1}$ be the image of the unit of $U_q(\mfm_t)$. Clearly $M_{\lambda,b}$ is a $\lambda$-basic module. 
\begin{Lem}\label{LemEigSpace}
The left ideal $I_{\lambda,b}$ is a proper left ideal, and 
\[
\msB = \{T_c^{+,(n)}\widehat{1}\mid n \geq 0\}\cup \{T_c^{-,(n)}\widehat{1}\mid n< 0\}
\]
forms a basis of $M_{\lambda,b}$. 
\end{Lem} 
\begin{proof}
The same technique as in the proof of Proposition \ref{PropBasicAdm} shows that $I_{\lambda,b}$ is spanned as a vector space by elements of the form 
\[
T_{c}^{+,(n)}P(A_b),\qquad T_c^{-,(n)}P(A_b),
\]
where $n \in \Z_{\geq 0}$ and $P$ a polynomial with $\lambda$ as a root. By Proposition \ref{PropActAc}, it follows that the $[c]$-eigenspace of $iB_t$ in $I_{\lambda,b}$ must be the space of such $P(A_b)$. Since $U_q(\mfm_t) = \mcO_q(S_t^2)\subseteq \mcO_q(SU(2))$ is a domain (i.e.\ has no zero divisors), the space of $P(A_b)$ does not contain the unit (as $A_b$ must generate a free polynomial algebra in one variable). 

Similarly, the domain property of $U_q(\mfm_t)$ implies that 
\[
T_{c}^{\pm,(n)} \notin \{T_{c}^{\pm,(n)}P(A_b)\mid P(\lambda)=0\},
\] 
so the $T_c^{\pm,(n)}\widehat{1}$ are all non-zero. As they are eigenvectors for $iB_t$ at different eigenvalues, they are linearly independent. 
\end{proof}

\begin{Def}
We call $M_{\lambda,b}$ the \emph{canonical $\lambda$-basic module centered at $[b]$}. 
\end{Def}

We write the basis of Lemma \ref{LemEigSpace} as
\[
e_{b+2n} = T_b^{+,(n)}\widehat{1},\qquad e_{b-2n} = T_b^{-,(n)}\widehat{1},\qquad n \in \Z_{\geq0},
\]
so that $e_c$ is an eigenvector for $iB_t$ at eigenvalue $[c]$. 

\begin{Prop}
Let $\lambda\in \R$. Then any $\lambda$-basic $U_q(\mfsl(2,\R)_t)$-module $V$ is a quotient of some $M_{\lambda,b}$. 
\end{Prop}
\begin{proof}
If $v \in V$ is a cyclic $iB_t$-eigenvector of $V$ at eigenvalue $[b]$, then it is clear from the universal construction of $M_{\lambda,b}$ that there exists a unique map 
\[
M_{\lambda,b} \rightarrow V
\]
which sends $\widehat{1}$ to $v$ and which is a  $U_q(\mfsl(2,\R)_t)$-module map.
\end{proof}

We now look at invariant forms on $M_{\lambda,b}$. 

\begin{Lem}
Let $\lambda,b \in \R$. There exists a unique $*$-invariant bilinear form on $M_{\lambda,b}$ satisfying $\langle e_b,e_b\rangle = 1$. 
\end{Lem}
The $*$-invariance means that we have a sesqui-linear form satisfying \eqref{EqStarInv}.
\begin{proof}
The uniqueness is clear, since the $e_{c}$ need to be orthogonal to each other for distinct $c$, and $e_b$ is a cyclic vector. 

To show existence, let $l_b$ be the linear functional on $M_{\lambda,b}$ determined by $l_b(e_c) = \delta_{b,c}$. Then we can define on $U_q(\mfm_t)$ as a $U_q(\mfm_t)$-module the $*$-invariant bilinear form 
\[
\langle x,y\rangle = l_b(x^*ye_b).
\]
We claim that this form is also $*$-invariant as a $U_q(\mfsl(2,\R)_t)$-module under the identification \eqref{EqConcId}. Indeed, let $\chi_b$ be the character on $U_q(\mfk_t)$ determined by $iB_t \mapsto [b]$. Then $\chi_b$ is a $*$-character. We then compute for $h\in U_q(\mfk_t)$ and $x,y\in U_q(\mfm_t)$ that 
\begin{eqnarray*} \langle x, hy\rangle &=& \tau(y_{(2)},h_{(1)}) l_b(x^*y_{(1)}h_{(2)}e_b) \\
&=& \chi_b(h_{(2)}) \tau(y_{(2)},h_{(1)}) l_b(x^*y_{(1)}e_b)\\ 
&=& \chi_b(h_{(2)}) \tau(S^{-1}(x_{(3)}^*)x_{(2)}^*y_{(2)},h_{(1)}) l_b(x_{(1)}^*y_{(1)}e_b)\\ 
&=& \chi_b(h_{(3)}) \tau(S^{-1}(x_{(3)}^*),h_{(1)})\tau(x_{(2)}^*y_{(2)},h_{(2)}) l_b(x_{(1)}^*y_{(1)}e_b)\\ 
&=& \tau(S^{-1}(x_{(2)}^*),h_{(1)}) l_b(h_{(2)}(x_{(1)}^*y)e_b)\\
&=& \chi_b(h_{(2)})\tau(S^{-1}(x_{(2)}^*),h_{(1)})l_b(x_{(1)}^*ye_b)\\
&=& \overline{\chi_b(h_{(2)}^*)}\overline{\tau(x_{(2)},h_{(1)}^*)}\langle x_{(1)},y\rangle \\
&=& \langle h^*x,y\rangle.
\end{eqnarray*}
Observe now that
\[
M_{\lambda,b} = U_q(\mfm_t)/I_{\lambda,b} = U_q(\mfm_t)/(\Omega_t - \lambda)U_q(\mfm_t),
\]
from which it follows directly that our sesqui-linear form descends to $M_{\lambda,b}$. 
\end{proof} 

Clearly the kernel of the $*$-invariant form on $M_{\lambda,b}$ is a submodule $N_{\lambda,b}$. We denote 
\begin{equation}\label{EqIrrQuot}
L_{\lambda,b} = M_{\lambda,b}/N_{\lambda,b},
\end{equation}
and we obtain an induced non-degenerate $*$-invariant form on $L_{\lambda,b}$. We determine when this form is positive definite, or, what is the same, when the form on $M_{\lambda,b}$ is positive semi-definite. We introduce the following terminology.

\begin{Def}
Let $\lambda\in \R$. We call a $\lambda$-basic $U_q(\mfsl(2,\R)_t)$-module $V$ \emph{unitary} if there exists a $*$-invariant positive definite form on $V$. We call $V$ a \emph{unitary basic} module if it is unitary $\lambda$-basic for some $\lambda$. 
\end{Def}

\begin{Prop}\label{PropIneq}
The canonical $(\lambda,b)$-basic module of $U_q(\mfsl(2,\R)_t)$ is unitary if and only if for all $n\geq0$
\[
0\leq \prod_{0\leq m \leq n}(\lambda + \qdif{a+b+2m+1} )(\qdif{a-b-2m-1}-\lambda),
\]
\[
0 \leq \prod_{0\leq m \leq n}(\lambda +\qdif{a+b-2m-1} )(\qdif{a-b+2m+1}-\lambda),
\]
where we write $t = [[a]] = q^a -q^{-a}$ (for a unique $a\in \R$). 
\end{Prop}
\begin{proof}
Denote 
\[
r_c = \langle e_c,e_c\rangle. 
\]
From \eqref{EqXYZT}, we see that 
\[
(T_c^+)^* = T_c^- = \frac{\qdif{c}}{\qdif{c+2}}T_{c+2}^- + C_1 A_{c+2} + C_2 T_{c+2}^+ + C_3
\]
for certain constants $C_i$. We deduce that 
\[
r_{c+2} = \langle e_c,(T_c^+)^*T_c^+e_c\rangle = \frac{\qdif{c}}{\qdif{c+2}}\langle e_c,T_{c+2}^-T_c^+e_c\rangle,\qquad c\in b+2\Z_{\geq0}. 
\]
By Lemma \ref{LemCommRelATT}, we can simplify the right hand side and obtain, recalling that $t= [[a]]$, 
\[
r_{c+2} =  \frac{\qdif{c}}{\qdif{c+2}}(\lambda + \qdif{a+c+1})(\qdif{a-c-1}-\lambda)r_c. 
\]
Similarly, $T_c^+ = \frac{\qdif{c}}{\qdif{c-2}}T_{c-2}^+$ plus other terms, and we obtain 
\[
r_{c-2} =  \frac{\qdif{c}}{\qdif{c-2}}(\lambda + \qdif{a+c-1})(\qdif{a-c+1}-\lambda)r_c,\qquad c\in b-2\Z_{\geq0}.
\]
\end{proof} 

It follows from the above proposition (and its proof) that for each $(\lambda,b)$ there is at most one unitary $(\lambda,b)$-basic module, which must then necessarily coincide with $L_{\lambda,b}$ as defined by \eqref{EqIrrQuot}. We also have the following.

\begin{Prop}
An $\wSL(2,\R)$-admissible $*$-representation is irreducible if and only if it is a unitary basic module. 
\end{Prop}
\begin{proof}
One direction follows from Proposition \ref{PropBasicAdm} and the comment just before the proposition. Conversely, if $V$ is a unitary basic module, then any $U_q(\mfsl(2,\R)_t)$-intertwiner must necessarily be trivial, as the $iB_t$-eigenspaces are one-dimensional and there is a cyclic $iB_t$-eigenvector. But if then $W$ is a non-zero submodule of $V$, it must be the direct sum of its intersections with the $iB_t$-eigenspaces. As we have a $*$-representation and as $V$ is a direct sum of its eigenspaces by the admissibility assumption, we can take the orthogonal projection of $V$ onto $W$, obtaining a non-zero intertwiner which is hence the identity map. This shows that $W=V$. 
\end{proof}

We now come to the main classification theorem.

\begin{Theorem}\label{TheoClass}
Any unitary basic module is unitarily equivalent to one of the following mutually inequivalent $*$-representations: 
\begin{enumerate}
\item The $(\lambda,b)$-basic module $L_{\lambda,b}$ for some $b \in (a-2,a]$ and 
\begin{equation}\label{EqIneqUnit}
 -\qdif{a+b'+1}<\lambda < \qdif{a-b-1},
\end{equation}
where $b' \in (-a-2,-a]$ is uniquely determined by $b'-b \in 2\Z$. In this case, the $e_{c}$ with $c\in b+2\Z$ form a basis of $L_{\lambda,b}$. 
\item The $(\lambda,b)$-basic module $L_{\lambda,b}$ with either $b>-a$ and $\lambda = - \qdif{a+b-1}$, or $b>a$ and $\lambda = \qdif{a-b+1}$. In this case, the $e_c$ with $c\in b+2\Z_{\geq0}$ form a basis of $L_{\lambda,b}$. 
\item The $(\lambda,b)$-basic module $L_{\lambda,b}$ with either $b<-a$ and $\lambda = -\qdif{a+b+1}$, or $b<a$ and $\lambda = \qdif{a-b-1}$. In this case, the $e_c$ with $c\in b-2\Z_{\geq0}$ form a basis of $L_{\lambda,b}$.
\item  The $(\lambda,b)$-basic module $L_{\lambda,b}$ with either $b=-a$ and $\lambda = -(q^{-1}+q)$, or $b=a$ and $\lambda = q^{-1}+q$. In this case, $L_{\lambda,b}$ is one-dimensional. 
\end{enumerate}
\end{Theorem}
\begin{proof}
Note first that, for $x,y\in \R$, 
\[
\qdif{x} \leq \qdif{y} \qquad \Leftrightarrow \qquad |x|\leq |y|.
\]
Secondly, it is clear from the proof of Proposition \ref{PropIneq} that if we denote by $C_{\lambda,b}$ the index values $c$ in $b+ 2\Z$ for which the basis vector $e_c \in M_{\lambda,b}$ has non-zero image in a unitarizable $L_{\lambda,b}$, then $C_{\lambda,b}$ must be of the form $b+ J_{\lambda,b}$ for some interval $J_{\lambda,b} \subseteq \Z$. 

If the indices have no upper or lower bound, any basis vector $e_b$ is cyclic, and hence $L_{\lambda,b} \cong L_{\lambda,b+2n}$ for any $n\in \Z$. We can hence assume that $b$ is such that $b\in (a-2,a]$, and the inequalities in Proposition \ref{PropIneq} give the desired inequality \eqref{EqIneqUnit}. 

If the indices have a lower bound but no upper bound, we are in the second case. Then we can choose our index $b$ such that $e_b \neq 0$ and $T_{b}^- e_b = 0$. This latter condition forces $\lambda = -\qdif{a+b-1}$ or $\lambda = \qdif{a-b+1}$. The other inequalities in Proposition \ref{PropIneq} then  lead to the respective conditions $b>-a$ and $b>a$. The third case is treated similarly.

Finally, if the indices have both a lower and upper bound, we must have either $\lambda =-\qdif{a+b+2m+1} = -\qdif{a+b-2n-1}$ for $m,n\in \Z_{\geq0}$ or $\lambda =\qdif{a-b-2m-1} = \qdif{a-b+2n+1}$ for $m,n\in \Z_{\geq 0}$. The remaining inequalities in Proposition \ref{PropIneq} are then easily seen to force $m=n$ and hence respectively $a=-b$ and $a=b$. 
\end{proof} 

In particular, we obtain the following. 

\begin{Theorem}\label{TheoMain}
Up to unitary equivalence, all the $SL(2,\R)_t$-admissible irreducible $*$-representations of $U_q(\mfsl(2,\R)_t)$ with $t =[[a]]$ are the following, all of them being mutually non-equivalent: 
\begin{enumerate}
\item The modules $\msL_{\lambda}^+ = L_{\lambda,a}$ with
\[
-(q^{1-2s} + q^{2s-1})<\lambda < q^{-1}+q,
\]
where $s \in [0,1)$ and $a+s\in \Z$.
\item The modules $\msL_{\lambda}^- = L_{\lambda,a-1}$ with 
\[
 -(q^{1-2s} + q^{2s-1})< \lambda < 2
\]
where $s\in [0,1)$ and $a+s\in \frac{1}{2}+\Z$.
\item The modules $\msD_n^+ = L_{\qdif{n-1},a+n}$ for $n \in \Z_{>0}$,
\item The modules $\msE_n^+ = L_{-\qdif{2a+n-1},a+n}$ for $n \in \Z_{>-2a}$,
\item The modules $\msD_n^- = L_{\qdif{n-1},a-n}$ for $n \in \Z_{>0}$,
\item The modules $\msE_n^- = L_{-\qdif{-2a+n-1},a-n}$ for $n \in \Z_{>2a}$,  
\item The trivial module $\msI = L_{q^{-1}+q,a}$.
\item The one-dimensional module $\msC = L_{-(q^{-1}+q),-a}$ if $a \in \frac{1}{2}\Z$.
\end{enumerate}
\end{Theorem}

Below, we represent the list of irreducible $SL(2,\R)_t$-admissible representations of quantum $SL(2,\R)_t$ graphically, see Figures \ref{Fig1}, \ref{Fig2}, \ref{Fig3}. The behaviour is $1$-periodic in the variable $a$ with $t = [[a]]$, and hence we can restrict to the values $a\in [0,1]$. The case for $1/2 < a <1$ is the same as for the case $0<a<1/2$, except that for the odd series the value $a$ is replaced by $1-a$. 

  \begin{center}    

\begin{figure}[hbtp!]
\centering
\begin{flushleft}
\underline{Even representations}
\end{flushleft}
\begin{tikzpicture}[scale = 0.7]

\filldraw[black, thick] (-5.4,0) -- (5.4,0);

\filldraw[black] (-5.5,0.3) circle (2pt) node[anchor=south]{$\msE_2^+$};
\filldraw[black] (-5.5,0) circle (2pt) node[anchor=east]{$\msC$} ;
\filldraw[black] (-5.5,-0.3) circle (2pt) node[anchor=north]{$\msE_2^-$};
\filldraw[black] (-6.3,1.1) circle (2pt) node[anchor=south]{$\msE_4^+$};
\filldraw[black] (-6.3,-1.1) circle (2pt) node[anchor=north]{$\msE_4^-$};

\draw (-4.2,0)  (-4.2,-0) node[anchor=south]{{\scriptsize $-(q+q^{-1})$}};
\draw (4.5,0)  (4.5,0) node[anchor=south]{{\scriptsize $q+q^{-1}$}};

\filldraw[black, dashed] (-6.8,1.6) -- (-7.7,2.5);
\filldraw[black, dashed] (-6.8,-1.6) -- (-7.7,-2.5);
\filldraw[black, dashed] (6.8,1.6) -- (7.7,2.5);
\filldraw[black, dashed] (6.8,-1.6) -- (7.7,-2.5);

\filldraw[black] (5.5,0.3) circle (2pt) node[anchor=south]{$\msD_2^+$};
\filldraw[black] (5.5,0) circle (2pt) node[anchor=west]{$\msI$};
\filldraw[black] (5.5,-0.3) circle (2pt) node[anchor=north]{$\msD_2^-$};
\filldraw[black] (6.3,1.1) circle (2pt) node[anchor=south]{$\msD_4^+$};
\filldraw[black] (6.3,-1.1) circle (2pt) node[anchor=north]{$\msD_4^-$};

\draw (-2.5,0.2) -- (-2.5,-0.2) node[anchor=north east]{{\scriptsize Complementary}};
\draw (-2.5,0.06)  (-2.5,0.06) node[anchor=south]{{\scriptsize $-2$}};
\draw (2.5,0.2) -- (2.5,-0.2) node[anchor=north west]{{\scriptsize {\scriptsize Complementary}}};
\draw (-2.5,0.06)  (2.5,0.06) node[anchor=south]{{\scriptsize $2$}};

\draw (0,0.2)  (0,-0.2)  node[anchor=north]{{\scriptsize Principal}};
\draw (0,0)  (0,0) node{$\times$};
\draw (0,0.22)  (0,0.22) node[anchor=south]{$0$};
\draw (0,0.9)  (0,0.9) node[anchor=south]{$\msL_{\lambda}^+$};

\end{tikzpicture}

\begin{flushleft}
\underline{Odd representations}
\end{flushleft}
\begin{tikzpicture}[scale = 0.7]

\filldraw[black, thick] (-3.4,0) -- (3.4,0);

\filldraw[black] (-3.5,0.15) circle (2pt) node[anchor=south]{$\msE_1^+$};
\filldraw[black] (-3.5,-0.15) circle (2pt) node[anchor=north]{$\msE_1^-$};
\filldraw[black] (-4.3,1) circle (2pt) node[anchor=south]{$\msE_3^+$};
\filldraw[black] (-4.3,-1) circle (2pt) node[anchor=north]{$\msE_3^-$};

\filldraw[black, dashed] (-4.8,1.5) -- (-6,2.5);
\filldraw[black, dashed] (-4.8,-1.5) -- (-6,-2.5);
\filldraw[black, dashed] (4.8,1.5) -- (6,2.5);
\filldraw[black, dashed] (4.8,-1.5) -- (6,-2.5);

\filldraw[black] (3.5,0.15) circle (2pt) node[anchor=south]{$\msD_1^+$};
\filldraw[black] (3.5,-0.15) circle (2pt) node[anchor=north]{$\msD_1^-$};
\filldraw[black] (4.3,1) circle (2pt) node[anchor=south]{$\msD_3^+$};
\filldraw[black] (4.3,-1) circle (2pt) node[anchor=north]{$\msD_3^-$};

\draw (0,0.2)  (0,-0.2)  node[anchor=north]{{\scriptsize Principal}};
\draw (0,0)  (0,0) node{$\times$};
\draw (0,0.22)  (0,0.22) node[anchor=south]{$0$};
\draw (-3,0.06)  (-3,0.06) node[anchor=south]{{\scriptsize {\scriptsize $-2$}}};
\draw (3,0.06)  (3,0.06) node[anchor=south]{{\scriptsize {\scriptsize $2$}}};
\draw (0,0.9)  (0,0.9) node[anchor=south]{$\msL_{\lambda}^-$};

\end{tikzpicture}
\caption{\text{Case $a=0$}}
\label{Fig1}
\end{figure}
\end{center}



\begin{figure}[hbtp!]
\centering
\begin{flushleft}
\underline{Even representations}
\end{flushleft}
\begin{tikzpicture}[scale = 0.7]

\filldraw[black, thick] (-3.4,0) -- (5.4,0);

\filldraw[black] (-3.5,0.15) circle (2pt) node[anchor=south]{$\msE_2^+$};
\filldraw[black] (-3.5,-0.15) circle (2pt) node[anchor=north]{$\msE_2^-$};
\filldraw[black] (-4.3,1) circle (2pt) node[anchor=south]{$\msE_4^+$};
\filldraw[black] (-4.3,-1) circle (2pt) node[anchor=north]{$\msE_4^-$};

\filldraw[black, dashed] (-4.8,1.5) -- (-5.7,2.5);
\filldraw[black, dashed] (-4.8,-1.5) -- (-5.7,-2.5);
\filldraw[black, dashed] (6.8,1.6) -- (7.7,2.5);
\filldraw[black, dashed] (6.8,-1.6) -- (7.7,-2.5);

\filldraw[black] (5.5,0) circle (2pt) node[anchor=west]{$\msI$};

\filldraw[black] (5.5,0.3) circle (2pt) node[anchor=south]{$\msD_2^+$};
\draw (4.4,0.2)  (4.4,0.2) node[anchor=south]{{\scriptsize $q+q^{-1}$}};
\filldraw[black] (5.5,-0.3) circle (2pt) node[anchor=north]{$\msD_2^-$};
\filldraw[black] (6.3,1.1) circle (2pt) node[anchor=south]{$\msD_4^+$};
\filldraw[black] (6.3,-1.1) circle (2pt) node[anchor=north]{$\msD_4^-$};

\draw (-3,0.06)  (-3,0.06) node[anchor=south]{{\scriptsize $-2$}};
\draw (2,0.2) -- (2,-0.2) node[anchor=north west]{{\scriptsize Complementary}};
\draw (-3,0.22)  (2,0.22) node[anchor=south]{{\scriptsize $2$}};
\draw (0,0.2)  (0,-0.2)  node[anchor=north]{{\scriptsize Principal}};
\draw (-0.5,0)  (-0.5,0) node{$\times$};
\draw (-0.5,0.22)  (-0.5,0.22) node[anchor=south]{$0$};
\draw (-0.5,0.9)  (-0.5,0.9) node[anchor=south]{$\msL_{\lambda}^+$};
\end{tikzpicture}


\begin{flushleft}
\underline{Odd representations}
\end{flushleft}
\begin{tikzpicture}[scale = 0.7]

\filldraw[black, thick] (-5.4,0) -- (3.4,0);

\filldraw[black] (-5.5,0.3) circle (2pt) node[anchor=south]{$\msE_1^+$};
\filldraw[black] (-5.5,-0.3) circle (2pt) node[anchor=north]{$\msE_1^-$};
\filldraw[black] (-6.3,1.1) circle (2pt) node[anchor=south]{$\msE_3^+$};
\filldraw[black] (-6.3,-1.1) circle (2pt) node[anchor=north]{$\msE_3^-$};

\filldraw[black] (-5.5,0) circle (2pt) node[anchor=east]{$\msC$};

\filldraw[black, dashed] (-6.8,1.6) -- (-7.7,2.5);
\filldraw[black, dashed] (-6.8,-1.6) -- (-7.7,-2.5);
\filldraw[black, dashed] (4.8,1.5) -- (5.7,2.5);
\filldraw[black, dashed] (4.8,-1.5) -- (5.7,-2.5);

\filldraw[black] (3.5,0.15) circle (2pt) node[anchor=south]{$\msD_1^+$};
\filldraw[black] (3.5,-0.15) circle (2pt) node[anchor=north]{$\msD_1^-$};
\filldraw[black] (4.3,1) circle (2pt) node[anchor=south]{$\msD_3^+$};
\filldraw[black] (4.3,-1) circle (2pt) node[anchor=north]{$\msD_3^-$};

\draw (-2,0.2) -- (-2,-0.2) node[anchor=north east]{{\scriptsize Complementary}};
\draw (-2,0.22)  (-2,0.22) node[anchor=south]{{\scriptsize $-2$}};
\draw (-2,0.06)  (3,0.06) node[anchor=south]{{\scriptsize $2$}};
\draw (-4.4,0.1)  (-4.4,0.1) node[anchor=south]{{\scriptsize$-(q+q^{-1})$}};
\draw (0.5,0.2)  (0.5,-0.2)  node[anchor=north]{{\scriptsize Principal}};
\draw (0.5,0)  (0.5,0) node{$\times$};
\draw (0.5,0.22)  (0.5,0.22) node[anchor=south]{$0$};
\draw (0.5,0.9)  (0.5,0.9) node[anchor=south]{$\msL_{\lambda}^-$};
\end{tikzpicture}
\caption{\text{Case $a$ half integer}}
\label{Fig2}

         

\vspace{0.8cm}
\begin{flushleft}
\underline{Even representations}
\end{flushleft}
\begin{tikzpicture}[scale = 0.7]

\filldraw[black, thick] (-5.4,0) -- (5.4,0);

\filldraw[black] (-5.5,0.15) circle (2pt) node[anchor=south]{$\msE_2^+$};
\draw (-3.6,0.1)  (-3.6,0.1) node[anchor= south]{{\scriptsize $-(q^{2a -1}+q^{-2a +1})$}};
\filldraw[black] (-5.5,-0.15) circle (2pt) node[anchor=north]{$\msE_2^-$};
\filldraw[black] (-6.3,1) circle (2pt) node[anchor=south]{$\msE_4^+$};
\filldraw[black] (-6.3,-1) circle (2pt) node[anchor=north]{$\msE_4^-$};
\filldraw[black] (5.5,0.3) circle (2pt) node[anchor=south]{$\msD_2^+$};
\filldraw[black] (5.5,0) circle (2pt) node[anchor=west]{$\msI$};

\filldraw[black, dashed] (-6.8,1.6) -- (-7.7,2.5);
\filldraw[black, dashed] (-6.8,-1.6) -- (-7.7,-2.5);
\filldraw[black, dashed] (6.8,1.4) -- (7.7,2.3);
\filldraw[black, dashed] (6.8,-1.4) -- (7.7,-2.3);

\filldraw[black] (5.5,-0.3) circle (2pt) node[anchor=north]{$\msD_2^-$};
\filldraw[black] (6.3,1.1) circle (2pt) node[anchor=south]{$\msD_4^+$};
\filldraw[black] (6.3,-1.1) circle (2pt) node[anchor=north]{$\msD_4^-$};

\draw (4.5,0.1)  (4.5,0.1) node[anchor= south]{{\scriptsize $q+q^{-1}$}};
\draw (-1.5,0.2) -- (-1.5,-0.2) node[anchor=north east]{{\scriptsize Complementary}};
\draw (-1.5,0.22)  (-1.5,0.22) node[anchor=south]{{\scriptsize $-2$}};
\draw (1.5,0.2) -- (1.5,-0.2) node[anchor=north west]{{\scriptsize Complementary}};
\draw (-1.5,0.22)  (1.5,0.22) node[anchor=south]{{\scriptsize $2$}};
\draw (0,0.2)  (0,-0.2)  node[anchor=north]{{\scriptsize Principal}};
\draw (0,0)  (0,0) node{$\times$};
\draw (0,0.22)  (0,0.22) node[anchor=south]{$0$};
\draw (0,0.9)  (0,0.9) node[anchor=south]{$\msL_{\lambda}^+$};
\end{tikzpicture}


\begin{flushleft}
\underline{Odd representations}
\end{flushleft}
\begin{tikzpicture}[scale = 0.7]
\filldraw[black, thick] (-5.4,0) -- (3.4,0);

\filldraw[black] (-5.5,0.15) circle (2pt) node[anchor=south]{$\msE_1^+$};
\filldraw[black] (-5.5,-0.15) circle (2pt) node[anchor=north]{$\msE_1^-$};
\filldraw[black] (-6.3,1) circle (2pt) node[anchor=south]{$\msE_3^+$};
\filldraw[black] (-6.3,-1) circle (2pt) node[anchor=north]{$\msE_3^-$};

\filldraw[black, dashed] (-6.8,1.6) -- (-7.7,2.5);
\filldraw[black, dashed] (-6.8,-1.6) -- (-7.7,-2.5);
\filldraw[black, dashed] (4.8,1.4) -- (5.7,2.2);
\filldraw[black, dashed] (4.8,-1.4) -- (5.7,-2.2);

\filldraw[black] (3.5,0.15) circle (2pt) node[anchor=south]{$\msD_1^+$};
\filldraw[black] (3.5,-0.15) circle (2pt) node[anchor=north]{$\msD_1^-$};
\filldraw[black] (4.3,1) circle (2pt) node[anchor=south]{$\msD_3^+$};
\filldraw[black] (4.3,-1) circle (2pt) node[anchor=north]{$\msD_3^-$};

\draw (-2,0.22)  (-2,0.22) node[anchor=south]{{\scriptsize $-2$}};
\draw (-2,0.22)  (3,0.22) node[anchor=south]{{\scriptsize $2$}};
\draw (-2,0.2) -- (-2,-0.2) node[anchor=north east]{{\scriptsize Complementary}};
\draw (-4.1,0.1)  (-4.1,0.1) node[anchor=south]{{\scriptsize $-(q^{2a}+q^{-2a})$}};
\draw (0.5,0.2)  (0.5,-0.2)  node[anchor=north]{{\scriptsize Principal}};
\draw (0.5,0)  (0.5,0) node{$\times$};
\draw (0.5,0.22)  (0.5,0.22) node[anchor=south]{$0$};
\draw (0.5,0.9)  (0.5,0.9) node[anchor=south]{$\msL_{\lambda}^-$};
\end{tikzpicture}
\caption{\text{Case $0 < a < \frac{1}{2}$ }}
\label{Fig3}
\end{figure}
\vfill

\newpage
\begin{Rem}
Let us determine in which way the above representations deform the usual irreducible representations of $SL(2,\R)$ \cite{Bar47}, pictorially represented in Figure \ref{Fig4}. For simplicity, we consider the case $t=0$. Then $B_0 = \widetilde{B}_0$, and a straightforward but tedious computation, using \eqref{EqDifEt}, shows that 
\begin{multline*}
(q-1)^{-2}(\Omega_0 - q-q^{-1}) = -2q^{-1}x+2z - 2(xz-yw) +q^{-1}(w^2+x^2) \\+q(z^2+q^{-2}y^2) -(q+1)(q^{-1/2}w+q^{-3/2}y)B_0 +  \mcO((q-1)).
\end{multline*}
Hence $\omega_q := (q-q^{-1})^{-2}(\Omega_0 - q-q^{-1})$ converges to the Casimir operator 
\[
\Omega = \frac{1}{4}H^2 -\frac{1}{2}H+EF
\] 
in $U(\mfsl(2,\R))$.

\begin{figure}[hbtp!]
\centering
\begin{flushleft}
\underline{Even representations}
\end{flushleft}
\begin{tikzpicture}[scale = 0.7]

\filldraw[black, thick] (-5,0) -- (5.5,0);
\filldraw[black] (5.5,0.3) circle (2pt) node[anchor=south]{$\msD_2^+$};
\filldraw[black] (5.5,0) circle (2pt) node[anchor=west]{$\msI$};
\filldraw[black, dashed] (-5.2,0) -- (-7,0);
\filldraw[black, dashed] (6.8,1.4) -- (7.7,2.3);
\filldraw[black, dashed] (6.8,-1.4) -- (7.7,-2.3);
\filldraw[black] (5.5,-0.3) circle (2pt) node[anchor=north]{$\msD_2^-$};
\filldraw[black] (6.3,1.1) circle (2pt) node[anchor=south]{$\msD_4^+$};
\filldraw[black] (6.3,-1.1) circle (2pt) node[anchor=north]{$\msD_4^-$};

\draw (1.7,0.2)  (1.7,-0.2) node[anchor=north west]{{\scriptsize Complementary}};
\draw (5.1,-0.6)  (5.1,-0.6) node[anchor=south]{{\scriptsize $0$}};
\draw (-2,0.2)  (-2,-0.2)  node[anchor=north]{{\scriptsize Principal}};
\draw (1,0)  (1,0) node{$\times$};
\draw (1,-0.1)  (1,-0.1) node[anchor=north]{$-1/4$};

\end{tikzpicture}


\begin{flushleft}
\underline{Odd representations}
\end{flushleft}
\begin{tikzpicture}[scale = 0.7]

\filldraw[black, thick] (-3,0) -- (5.5,0);
\filldraw[black] (5.5,0.3) circle (2pt) node[anchor=south]{$\msD_1^+$};
\filldraw[black, dashed] (-3.2,0) -- (-5,0);
\filldraw[black, dashed] (6.8,1.4) -- (7.7,2.3);
\filldraw[black, dashed] (6.8,-1.4) -- (7.7,-2.3);
\filldraw[black] (5.5,-0.3) circle (2pt) node[anchor=north]{$\msD_1^-$};
\filldraw[black] (6.3,1.1) circle (2pt) node[anchor=south]{$\msD_3^+$};
\filldraw[black] (6.3,-1.1) circle (2pt) node[anchor=north]{$\msD_3^-$};

\draw (4.8,-0.7)  (4.8,-0.7) node[anchor=south]{{\scriptsize $-1/4$}};
\draw (0.5,0.2)  (0.5,-0.2)  node[anchor=north]{{\scriptsize Principal}};

\end{tikzpicture}
\caption{\text{Case of classical $SL(2,\R)$} (parametrisation by values of $\Omega$)}
\label{Fig4}
\end{figure}

Now our condition for the first case in Theorem \ref{TheoClass} can be rewritten 
\[
-\frac{\qdif{\frac{b}{2}}\qdif{\frac{b}{2}+1}}{(q-q^{-1})^2} < \omega_q < \left[\frac{b}{2}\right]\left[\frac{b}{2}+1\right]. 
\]
We see that the left hand side condition goes to minus infinity as $q$ tends to $1$, and in particular in that limit the representations $\msE_n^{\pm}$ and $\msC$ disappear. A similar phemenon is known to happen in the setting of quantum $SL(2,\C)$ \cite{Pus93,PW94}. Of the remaining representations, the $\msD_{n}^{\pm}$ with $n\geq 2$ are the analogues of the positive and negative discrete series, with $\msD_1^{\pm}$ corresponding to the limiting positive and negative discrete series representations. The $\msL^{\pm}_{\lambda}$ are analogues of the positive and negative principal series for $\lambda \leq 2$, with a complementary series given as $\msL_{\lambda}^+$ for $2<\lambda < q^{-1}+q$. It will be motivated in future work that also the $\msL_{\lambda}^{\pm}$ with $\lambda< -2$ should also be interpreted as complementary series representations in the quantum setting (vanishing as $q$ tends to $1$).
\end{Rem}

\subsection{Admissible representations}

In this section, we briefly return to the setting of general admissible modules, and clarify the link with admissible representations. First, we construct the following $\widetilde{SL}(2,\R)$-admissible modules. 

\begin{Lem}\label{LemNewForm}
Let $\lambda\in \R$ and $\chi \in \R/2\Z$, where we concretely view $\R/2\Z$ as coset classes in $\R$. Then there exists a vector space $V_{\lambda,\chi}^+$ with basis $\{\xi_{c}\mid c\in\chi\}$ carrying a $\widetilde{SL}(2,\R)$-admissible $U_q(\mfsl(2,\R)_t)$-module structure such that each $\xi_c$ is an $iB_t$-eigenvector at eigenvalue $[c]$ and such that 
\[
A_c\xi_c = \lambda \xi_c,\qquad T_c^{\pm}\xi_c = (\qdif{a\pm c+1}\pm\lambda)\xi_{c\pm 2}. 
\]
Similarly, there exists a vector space $V_{\lambda,\chi}^-$ with basis $\{\xi_{c}\mid c\in \chi\}$ carrying a $\widetilde{SL}(2,\R)$-admissible $U_q(\mfsl(2,\R)_t)$-module structure such that each $\xi_c$ is an $iB_t$-eigenvector at eigenvalue $[c]$ and such that 
\[
A_c\xi_c = \lambda \xi_c,\qquad T_c^{\pm}\xi_c = (\qdif{a\mp c-1}\mp\lambda)\xi_{c\pm 2}. 
\]
\end{Lem} 
\begin{proof}
Assume first that $\lambda\geq0$, and choose $b\in \chi$ such that $\qdif{a-c+1}\neq \lambda$ for all $c\in b-2\Z_{\geq0}$. Consider inside $M_{\lambda,b}$ the following vectors for $n\in \Z_{\geq 0}$, 
\begin{multline*}
\xi_{b+2n} = \left(\prod_{m=0}^{n-1} (\qdif{a+b+2m+1}+\lambda)\right)^{-1}e_{b+2n},\\ \xi_{b-2n} = \left(\prod_{m=0}^{n-1}(\qdif{a-b+2m-1}-\lambda)\right)^{-1}e_{b-2n}. 
\end{multline*}
Then one easily checks that this basis provides a realisation of $V_{\lambda,\chi}^+$. A similar rescaling can be done in the case $\lambda <0$, or to obtain the representation $V_{\lambda,\chi}^-$.  
\end{proof}

\begin{Prop}\label{PropAdmRep}
Providing $V_{\lambda,\chi}^{\pm}$ with the pre-Hilbert space structure for which the $\xi_{c}$ from Lemma \ref{LemNewForm} form an orthonormal basis, it becomes an $\widetilde{SL}(2,\R)$-admissible representation.  
\end{Prop} 

\begin{proof}
Using \eqref{EqXYZT}, we obtain that the action of the generators $X_t,Z_t,Y_t$ on $V_{\lambda,\chi}^{\pm}$ is given concretely by
\begin{multline*}
iX_t\xi_c = q^{c+1}\frac{\qdif{a\pm (c+1)}\pm\lambda}{\qdif{c}\qdif{c+1}}\xi_{c+2}+ \frac{[[c]][[a]]+\qdif{2}\lambda}{\qdif{c-1}\qdif{c+1}}\xi_c \\ - q^{-c+1}\frac{\qdif{a\mp (c-1)}\mp \lambda}{\qdif{c-1}\qdif{c}}\xi_{c-2}. 
\end{multline*}
\begin{multline*}
Z_t\xi_c = -\frac{\qdif{a\pm(c+1)}\pm \lambda}{\qdif{c}\qdif{c+1}}\xi_{c+2} + \frac{-\qdif{2}[[a]] +[[c]]\lambda}{\qdif{c-1}\qdif{c+1}}\xi_c \\ - \frac{\qdif{a\mp(c-1)}\mp \lambda}{\qdif{c-1}\qdif{c}}\xi_{c-2}, 
\end{multline*}
\begin{multline*}
iY_t\xi_c = q^{-c-1} \frac{\qdif{a\pm (c+1)}\pm \lambda}{\qdif{c}\qdif{c+1}}\xi_{c+2} - \frac{[[c]][[a]]+\qdif{2}\lambda}{\qdif{c-1}\qdif{c+1}}\xi_c \\-q^{c-1} \frac{\qdif{a\mp(c-1)}\mp \lambda}{\qdif{c-1}\qdif{c}}\xi_{c-2}. 
\end{multline*} 
As the coefficients of $X_t,Z_t,Y_t$ are uniformly bounded when $c$ varies, it follows that they define bounded operators. 
\end{proof} 

\begin{Cor}
Any irreducible $\widetilde{SL}(2,\R)$-admissible module carries a pre-Hilbert space structure for which it becomes an irreducible $\widetilde{SL}(2,\R)$-admissible representation. 
\end{Cor} 

\begin{proof}
Since any irreducible $\widetilde{SL}(2,\R)$-admissible module is a quotient of some $M_{\lambda,b}$, it is sufficient to show that the $M_{\lambda,b}$ can be made into $\widetilde{SL}(2,\R)$-admissible representations. We give the proof for $\lambda\geq 0$, the proof for $\lambda <0$ is similar. 

By the proof of Lemma \ref{LemNewForm}, we see that $M_{\lambda,b}\cong V_{\lambda,b+2\Z}^+$ for $b\ll 0$, which can be made into an  irreducible $\widetilde{SL}(2,\R)$-admissible representation by Proposition \ref{PropAdmRep}. To prove that the corollary holds for an arbitrary $M_{\lambda,b}$, we can hence by induction assume that it already holds for all $M_{\lambda,c}$ with $c\in b+2\Z$ and $c<b$. 

Consider the natural intertwiner $M_{\lambda,b}\overset{\phi}{\rightarrow} V_{\lambda,b+2\Z}^+$ with $e_b \mapsto \xi_b$. As $\lambda \geq 0$, we see by the specific form of the action of the $T_{c}^+$ on $V_{\lambda,b+2\Z}^+$ that $\phi$ is surjective. There hence exists a largest $c\in b+2\Z$ with $c<b$ such that $e_c\in \Ker(\phi)$. Consider the natural intertwiner $M_{\lambda,c}\overset{\psi}{\rightarrow} M_{\lambda,b}$ with $e_c \mapsto e_c$. Then this map must land in $\Ker(\phi)$, and must be surjective by considering the action of the $T_{c}^{-}$. By induction, $M_{\lambda,c}$ can be made into an irreducible $\widetilde{SL}(2,\R)$-admissible representation, and hence the same holds for $\Ker(\phi)$. As this also holds for $V_{\lambda,b+2\Z}^+$, the same must hold for $M_{\lambda,b}$. 
\end{proof}

It is also not hard to classify the finite-dimensional irreducible $\widetilde{SL}(2,\R)_t$-admissible $U_q(\mfsl(2,\R)_t)$-modules. 

\begin{Prop}
For each $N\in \Z_{>0}$, there exist exactly two finite-dimensional irreducible $\widetilde{SL}(2,\R)_t$-admissible $U_q(\mfsl(2,\R)_t)$-modules $V_N^+$ and $V_N^-$ of dimension $N$, which are $\lambda$-basic for respectively $\lambda = \qdif{N}$ and $\lambda = -\qdif{N}$. The modules $V^+_N$ are automatically $SL(2,\R)_t$-admissible, while for any $N$ the module $V^-_N$ is $SL(2,\R)_t$-admissible if and only if $a$ is a half-integer. 
\end{Prop} 
\begin{proof}
Let $(V,\pi)$ be a finite-dimensional irreducible $\widetilde{SL}(2,\R)_t$-admissible $U_q(\mfsl(2,\R)_t)$-module, which we may assume to be $\lambda$-basic for some $\lambda \in \R$. Necessarily, by irreducibility and the way the $T_c^{\pm}$ act on $iB_t$-eigenvectors, the eigenvalues of $\pi(iB_t)$ must be of the form $\{[b],[b+2],\ldots,[b+2N-2]\}$ for some $b\in \R$. 

We hence have a unique surjective intertwiner $M_{\lambda,b}\rightarrow V$ with $e_b \mapsto \xi_b$, for some non-zero eigenvector $\xi_b\in V$ at eigenvalue $[b]$. As $e_{b-2}$ and $e_{b+2N}$ lie in the kernel of this map, we must have $T_{b-2}^+ e_{b-2} = 0$ and $T_{b+2N}^- e_{b+2N} = 0$. By Lemma \ref{LemCommRelATT}, we deduce that 
\begin{multline*}
0 = (\lambda -  \qdif{a-b+1})(\lambda+ \qdif{a+b+1}) \\ = (\lambda - \qdif{a-b-2N+1})(\lambda + \qdif{a+b+2N+1}),
\end{multline*}
so either 
\[
\lambda = \qdif{a-b+1} = \qdif{a-b-2N+1},
\]
hence $b = a-N+1$ and $\lambda = \qdif{N}$, or 
\[
\lambda =  -\qdif{a+b+1} = -\qdif{a+b+2N+1},
\]
hence the $b = -a-N-1$ and $\lambda = -\qdif{N}$. Conversely, it is clear using \eqref{EqXYZT} that for $b,\lambda$ of the above form, the subspace of $M_{\lambda,b}$ spanned by the $e_c$ with $c<b$ or $c\geq b+2N$ is indeed a submodule. 

This proves the existence and uniqueness of the $N$-dimensional modules $V_N^{\pm}$. Clearly the $V_{N}^+$ are $SL(2,\R)_t$-admissible, as the eigenvalues of $iB_t$ are of the form $a+ M$ for an integer $M$. In the case of $V_N^-$ the eigenvalues  of $iB_t$ are of the form $-a+M$ for $M$ integer, hence $V_N^-$ is $SL(2,\R)_t$-admissible if and only if $a$ and $-a$ differ by an integer, i.e.~ when $a$ is a half-integer. 
\end{proof}

To end, let us explicitly relate the modules $V_{\lambda,b+2\Z}^{\pm}$ to the irreducible unitary representations obtained in Theorem \ref{TheoClass}. The proof is straightforward.

\begin{Prop}
Any irreducible $\widetilde{SL}(2,\R)_t$-admissible $*$-representation can be realized as a subquotient of some $V_{\lambda,\chi}^{\pm}$ with $\chi \in \R/2\Z$. More precisely, denote $\lambda_{\chi}^{\pm}  = \min_{c\in \chi}\{ \qdif{a \mp (c+1)}\}$.  
\begin{enumerate}
\item If $-\lambda_{\chi}^{-} <  \lambda < \lambda_{\chi}^{+}$, then $\lambda_{\chi}^+ = \qdif{a-b-1}$ for a unique $b\in (a-2,a]$, and $V_{\lambda,\chi}^+ \cong V_{\lambda,\chi}^- \cong L_{\lambda,b}^+$. 
\item If there exists $b\in \chi$ with $\lambda = \mp \qdif{a\pm (b-1)}$ and $b>\mp a$, then $L_{\lambda,b} \subseteq V_{\lambda,\chi}^{\pm}$ as the submodule generated by the $\xi_{b+2n}$ with $n \in \Z_{\geq0}$. 
\item If there exists $b\in \chi$ with $\lambda = \pm \qdif{a\mp (b+1)}$ and $b< \pm a$, then $L_{\lambda,b} \subseteq V_{\lambda,\chi}^{\pm}$ as the submodule generated by the $\xi_{b-2n}$ with $n \in \Z_{\geq0}$. 
\item If $\lambda = q^{-1}+q$ and $a \in\chi$, then $V_{\lambda,\chi}^{\pm}$ has a submodule spanned by the $\{\xi_{a \mp 2n}\mid n\in \Z_{\geq0}\}$, and 
\[
L_{\lambda, a} \cong \Span\{\xi_{a\mp 2n}\mid n\in \Z_{\geq0}\}/\Span\{\xi_{\pm a\mp 2n}\mid n\in \Z_{>0}\}.
\]
Similarly, if $\lambda = -(q^{-1}+q)$ and $-a \in\chi$, then $V_{\lambda,\chi}^{\pm}$ has a submodule spanned by the $\{\xi_{-a \pm 2n}\mid n\in \Z_{\geq0}\}$, and 
\[
L_{\lambda, a} \cong \Span\{\xi_{-a\pm 2n}\mid n\in \Z_{\geq0}\}/\Span\{\xi_{-a\pm 2n}\mid n\in \Z_{>0}\}.
\]
\end{enumerate}
 \end{Prop} 

\section{Outlook and perspectives}

In this paper, we have looked at a pair of Hopf $*$-algebras $A,U$ in duality, and have shown that if $B \subseteq A$ and $I \subseteq U$ are resp.\ a right and a left coideal $*$-subalgebra which are orthogonal to each other, they generate within the Drinfeld double Hopf $*$-algebra $\mcD(A,U)$ a right coideal $*$-subalgebra $\mcD(B,I)$, which we dub the  \emph{Drinfeld double coideal}. 

When considering the case $U = U_q(\mfsu(2))$ and $A = \mcO_q(SU(2))$, we have shown that if we take $I$ the coideal $*$-subalgebra $U_q(\mfk_t)$ generated by a particular element $B_t \in U_q(\mfsu(2))$ and $B= \mcO_q(T_t\backslash SU(2))$ the orthogonal coideal in $\mcO_q(SU(2))$, the resulting Drinfeld double coideal $\mcD(B,I) = U_q(\mfsl(2,\R)_t)$ can be seen as a quantization of the Hopf $*$-algebra $U(\mfsl(2,\R))$. The parameter $t$ encodes the particular position of (the Cartan subalgebra of) $\mfsl(2,\R)$ inside $\mfsl(2,\C)$. We have then classified all irreducible $*$-representations of $U_q(\mfsl(2,\R)_t)$, and more generally all admissible simple modules of $U_q(\mfsl(2,\R)_t)$.  

There are several directions in which this research can be continued. 
\begin{itemize}
\item The construction of Drinfeld double coideals can be considered in the C$^*$-algebraic and von Neumann algebraic framework of locally compact quantum groups \cite{KV00,KV03}. In the particular setting of discrete/compact quantum groups, this research has already been initiated in \cite{DCDz21}, where we also look at a theory of associated \emph{quasi-invariant integrals}.   
\item As indicated in the introduction, the procedure in this paper can be applied to construct quantizations of \emph{any} real semisimple Lie algebra. Indeed, if $\mfl$ is a real semisimple Lie algebra, we can consider its complexification $\mfg = \mfl \underset{\R}{\otimes} \C$. We can choose for $\mfg$ a compact real form  $\mfu \subseteq \mfg$ which is compatible with the real form $\mfl$, in the sense that $\mfk = \mfu \cap \mfl$ has maximal dimension. We can then choose a maximal Cartan subalgebra $\mft$ of $\mfu$ such that $\mft \cap \mfk$ has \emph{minimal dimension}. In the works of G.\ Letzter \cite{Let99,Let02} (see also the work of S.\ Kolb who treated the non-finite case \cite{Kol14}), it was shown how \emph{any} such pair can be quantized as $U_q(\mfk) \subseteq U_q(\mfu)$ with $U_q(\mfk)$ a coideal subalgebra. In \cite{DCNTY19} (see also \cite{DCM20}), it was shown how by a small modification we can make these into coideal $*$-subalgebras. With $U$ the simply connected Lie group integrating $\mfu$, denote by $\mcO_q(K\backslash U)$ the coideal $*$-subalgebra of $\mcO_q(U)$ orthogonal to $U_q(\mfk)$. By the same reasoning as in this paper, Drinfeld duality will then give that the resulting Drinfeld coideal $*$-subalgebra $U_q(\mfl) := \mcD(\mcO_q(K\backslash U),U_q(\mfk))$ will give a quantization of the universal enveloping $*$-algebra $U(\mfl)$ of $\mfl$. Of course, the study of the representation theory of $U_q(\mfl)$ will need the development of more sophisticated tools than the rather ad hoc ones used in this paper. 
\item Although they are constructed by very different methods, the representation theory of $U_q(\mfsl(2,\R))$ is in the end very similar to the (analytic) representation theory of extended $U_q(\mfsu(1,1))$, as treated in \cite{KK03}. Understanding the connection between these two representation theories is a very interesting problem, that should also shed light on the quantization problem for real semisimple Lie groups as locally compact quantum groups.  
\item Given that the representation theory of quantized $U_q(\mfsl(2,\C))$ (as a Drinfeld double of $U_q(\mfsu(2))$) is also explicitly known \cite{Pus93}, we can consider the associated \emph{branching problem} from representations of $U_q(\mfsl(2,\C))$ to representations of $U_q(\mfsl(2,\R)_t)$, using perhaps the more analytic machinery developed in \cite{DCDz21}. The latter paper also opens up the possibility to study the \emph{monoidal structure} of the representation theory of $U_q(\mfsl(2,\R)_t)$. 
\item We have not touched explicitly on induction techniques in our framework. Indeed, our classification basically hinges on an ad hoc induction from the component $U_q(\mfk_t)$ of $U_q(\mfsl(2,\R)_t)$, corresponding to induction from a maximal compact subgroup, but leaves open the more intrinsic construction of e.g.\ principal series representations using \emph{parabolic} induction. It is at the moment not clear how the latter can be achieved in full generality for quantizations of real Lie groups of arbitrary rank. For the case of quantum $SL(2,\R)$ this will be explored in future work, where also the quantum analogue of the Plancherel formula will be considered.   
\end{itemize}


\end{document}